\documentclass[]{article}

\usepackage{vmargin}
\setmarginsrb{2cm}{1.5cm}{2cm}{4cm}%
              {10pt}{10pt}{15pt}{21pt}          

\usepackage{amsmath,amssymb,amsfonts}
\usepackage{graphicx,epsfig,color,psfrag,xspace}
\usepackage[noend]{algorithmic}

\newcommand{\todo}[1]{\par\noindent{\color{red}\raggedright\textsc{#1}\par\marginpar{\Large$\star$}}}
\newcommand{\margin}[1]{\marginpar{\tiny\color{red} #1}}

\newcommand{\until}[1]{\{1,\dots, #1\}}
\newcommand{\subscr}[2]{#1_{\textup{#2}}}

\newcommand{\setdef}[2]{\{#1 \, | \; #2\}}
\newcommand{\seqdef}[2]{\{#1\}_{#2}}
\newcommand{\bigsetdef}[2]{\big\{#1 \, | \; #2\big\}}
\newcommand{\map}[3]{#1: #2 \rightarrow #3}
\newcommand{\setmap}[3]{#1: #2 \rightrightarrows #3}
\newcommand{\union}{\operatorname{\cup}}
\newcommand{\intersection}{\ensuremath{\operatorname{\cap}}}
\newcommand{\intersect}{\ensuremath{\operatorname{\cap}}}

\newcommand{\ON}{\textup{\texttt{ON}}\xspace}
\newcommand{\OFF}{\textup{\texttt{OFF}}\xspace}

\newcommand{\norminf}[1]{\ensuremath{\|#1\|_{\infty}}}

\newcommand\oprocendsymbol{\hbox{$\square$}}
\newcommand\oprocend{\relax\ifmmode\else\unskip\hfill\fi\oprocendsymbol}

\def\C{\mathcal{C}}
\def\CL{\mathcal{C}_{(\ell)}}

\newcommand\Prob{\mathbb{P}}
\newcommand\sat{\operatorname{sat}}
\def\part{\mathcal{V}_N}
\newcommand{\Cd}{\operatorname{Cd}}      
\newcommand{\CC}{\operatorname{CC}}      

\newcommand{\argmin}{\mathop{\operatorname{argmin}}}
\newcommand{\Tdelta}{\ensuremath{T^\delta}}
\newcommand{\integersnonnegative}{\ensuremath{\mathbb{Z}}_{\ge 0}}
\newcommand{\integernonnegative}{\ensuremath{\mathbb{Z}}_{\ge 0}}
\newcommand{\real}{\ensuremath{\mathbb{R}}}
\newcommand{\realpositive}{\ensuremath{\mathbb{R}}_{>0}}
\newcommand{\realnonnegative}{\ensuremath{\mathbb{R}}_{\geq0}}
\renewcommand{\natural}{{\mathbb{N}}}

\newcommand{\halfs}{\subscr{H}{bisector}}

\newcommand{\Hgeneric}{\subscr{\mathcal{H}}{multicenter}}
\newcommand{\Hvoronoi}{\subscr{\mathcal{H}}{Voronoi}}
\newcommand{\Hcentroid}{\subscr{\mathcal{H}}{centroid}}
\newcommand{\Hone}{\subscr{\mathcal{H}}{1}}
\newcommand{\HHdc}{\subscr{\mathcal{H}}{worst}}

\newcommand{\diam}{\operatorname{diam}}             
\newcommand{\dist}{\operatorname{dist}}             
\newcommand{\pdist}{\operatorname{pseudodist}}             
\newcommand{\interior}{\operatorname{int}}      
\newcommand{\closure}[1]{\overline{#1}}              
\newcommand{\measure}{\mu}          
\newcommand{\eps}{\epsilon}

\newcommand{\myset}{\setdef{(i,j)\in\until{N}^2}{i\not=j}}
\newcommand{\rcomm}{\subscr{r}{comm}}
\newcommand{\lcomm}{\subscr{\lambda}{comm}}
\newcommand{\pcomm}{\subscr{p}{comm}}
\newcommand{\parea}{\subscr{p}{area}}
\newcommand{\duration}{d}
\newcommand{\ee}{\textup{e}}

\newcommand{\Hequal}{\stackrel{\text{{\tiny(H)}}}{=}}

\newtheorem{theorem}{Theorem}[section]
\newtheorem{proposition}[theorem]{Proposition}
\newtheorem{corollary}[theorem]{Corollary}
\newtheorem{definition}[theorem]{Definition}
\newtheorem{lemma}[theorem]{Lemma}

\newtheorem{remark}[theorem]{Remark}
\newtheorem{remarks}[theorem]{Remarks}
\newtheorem{example}[theorem]{Example}

\newcommand{\qed}{\hfill \mbox{\raggedright \rule{.07in}{.1in}}}
\newenvironment{proof}{\vspace{1ex}\noindent{\it Proof.}\hspace{0.5em}}
	{\hfill\qed\vspace{1ex}}

        \title{Gossip Coverage Control for Robotic Networks: \\ Dynamical
          Systems on the Space of Partitions\thanks{ This work was
            supported in part NSF grant IIS-0904501, ARO MURI grant
            W911NF-05-1-0219, and ONR grant N00014-07-1-0721 A preliminary
            and incomplete version of this work appeared in the Proceedings
            of the 2009 American Control Conference, pages 2228-2235, St.\
            Louis, Missouri, USA. The authors gratefully acknowledge Prof.\
            Andrew R.\ Teel for Remark~\ref{remark:AndyTeel}. }}

\author{Francesco Bullo$^\dagger$\and Ruggero Carli\thanks{Francesco Bullo
    and Ruggero Carli are with the Center for Control, Dynamical Systems
    and Computation, University of California at Santa Barbara, Santa
    Barbara, CA 93106, USA,
    {\tt\small\{bullo,carlirug\}@engineering.ucsb.edu}.}\and Paolo
  Frasca\thanks{Paolo Frasca was with the Dipartimento di Matematica,
    Politecnico di Torino, and is now with the {D.I.I.M.A.}, Universit\`{a}
    di Salerno and with the Istituto per le Applicazioni del
    Calcolo\textendash{}{C.N.R.}, Roma, Italy, {\tt\small
      paolo.frasca@gmail.com}.}  }

\begin{document}
\maketitle

\begin{abstract}
  Future applications in environmental monitoring, delivery of services and
  transportation of goods motivate the study of deployment and partitioning
  tasks for groups of autonomous mobile agents. These tasks are achieved by
  recent coverage algorithms, based upon the classic methods by
  Lloyd. These algorithms however rely upon critical requirements on the
  communication network: information is exchanged synchronously among all
  agents and long-range communication is sometimes required.  This work
  proposes novel coverage algorithms that require only \emph{gossip
    communication}, i.e., asynchronous, pairwise, and possibly unreliable
  communication.  Which robot pair communicates at any given time may be
  selected deterministically or randomly.  A key innovative idea is
  describing coverage algorithms for robot deployment and environment
  partitioning as dynamical systems on a space of partitions.  In other
  words, we study the evolution of the regions assigned to each agent
  rather than the evolution of the agents' positions.  The proposed gossip
  algorithms are shown to converge to centroidal Voronoi partitions under
  mild technical conditions.

  Our treatment features a broad variety of results in topology, analysis
  and geometry. First, we establish the compactness of a suitable space of
  partitions with respect to the symmetric difference metric. Second, with
  respect to this metric, we establish the continuity of various geometric
  maps, including the Voronoi diagram as a function of its generators, the
  location of a centroid as a function of a set, and the widely-known
  multicenter function studied in facility location problems. Third, we
  prove two convergence theorems for dynamical systems on metric spaces
  described by deterministic and stochastic switches.
\end{abstract}

\thispagestyle{plain}
\section{Introduction}\label{sec:Intro}
In the not too distant future, networks of coordinated autonomous robots
will perform a broad range of environmental monitoring and logistic tasks.
Robotic camera networks will monitor airports and other public
infrastructures. Teams of vehicles will perform surveillance, exploration
and search and rescue operations. Groups of robots will enable novel
logistic capacities in the transportation of goods and the delivery of
services and resources to customers. New applications will be enabled by
the ongoing decreases in size and cost and the increases in performance of
sensors, actuators, communication devices and computing elements.

In these future applications, load balancing algorithms will dictate how
the workload is shared amongst and assigned to the individual robots. In
other words, robotic resources will be assigned and deployed to competing
requests in such a way as to optimize some performance metric. Remarkably,
load balancing problems in robotic networks are often equivalent to robotic
deployment and environment partitioning problems. For example, in
surveillance applications, optimal sensor coverage is often achieved by
partitioning the environment and assigning individual robotic sensors to
individual regions of responsibility.  Similarly, in the transportation of
goods or delivery of services, minimizing the customer wait-time is
equivalent to a multi-vehicle routing problem and, in turn, to computing
optimal depot positions and regions of responsibility.

Motivated by these scenarios, this paper considers the two following
interrelated problems.  The \emph{deployment problem} for a robotic network
amounts to the design of coordination algorithms that lead the robots to be
optimally placed in an environment of interest. Deployment performance is
characterized by an appropriate network utility function that measures the
deployment quality of a given configuration.  The \emph{partitioning
  problem} is the design of coordination algorithms that lead the robots to
optimally partition the environment into subregions of interest; even here
the objective is usually achieved through the design of appropriate utility
functions.

\subsection*{Literature review} A wide range of literature is relevant to
our investigation. In next four paragraphs we review previous work on
centroidal Voronoi partitions, distributed robotic deployment, animal
territory partitioning, and mathematical modeling tools.

The ``centering and partitioning'' algorithm originally proposed by
Lloyd~\cite{SPL:82} and elegantly reviewed in the survey~\cite{QD-VF-MG:99}
is a classic approach to facility location and environment partitioning
problems.  The Lloyd algorithm computes centroidal Voronoi partitions as
optimal configurations of an important class of objective functions, called
multicenter functions. Besides their intended application to quantization
theory~\cite{RMG-DLN:98}, centroidal Voronoi partitions have widespread
applications in numerous disciplines, including statistical pattern
recognition~\cite{AKJ-RPWD-JM:00}, mathematical
imaging~\cite{QD-MG-LJ-XW:06}, geometric optimization~\cite{PKA-MS:98} and
spatial resource allocation~\cite{ZD-HWH:01}.  Recent mathematical interest
has focused on convergence analysis~\cite{QD-ME-JL:06}, bifurcation
analysis of low dimensional problems~\cite{RV-DAWB-TG-TH-MF-JM:09}, and
anisotropic partitions~\cite{QD-DW:05}, among other topics.

Distributed and robotic versions of the Lloyd algorithm have been recently
developed in the multiagent literature; see the survey~\cite{SM-JC-FB:04n}
and the text~\cite[Chapter 5 and literature notes in Section
5.4]{FB-JC-SM:09}. We briefly review this growing literature in what
follows.
Generalized centroidal Voronoi partitions are shown in~\cite{RG-JC:07} to
be asymptotically optimal for estimation of stochastic spatial fields by
sensor networks.
In~\cite{JC-FB:08r} convergence results are obtained for simple basic
interactions among mobile agents such as ``move away from the closest other
agent'' or ``move toward the furthest vertex of your region of
responsibility.''
Convergence to centroidal Voronoi partitions is established
in~\cite{AA-KS-EF:09} for a class of communication-less sensor-based
algorithms (related to the classic clustering work by
MacQueen~\cite{JMC:67}). %
In~\cite{MS-DR-JJS:08} adaptive coverage controls are proposed for
environments described by unknown density functions.  
Dynamic environments and corresponding dynamic coverage problems are
treated in~\cite{IIH-DMS:07}. 
In~\cite{MP-EF-FB:08h} 
partitioning policies are shown to achieve optimal load balancing in
vehicle routing problems, i.e., problems in which a robotic network
provides service to customers that arrive in real time in the environment.

Territory partitioning via competitive behaviors is a classic subject of
study in behavioral ecology; see~\cite{ESA:01} for a comprehensive
survey. For example, it is known~\cite{FRA-DMG:03} that the foraging
behavior of conflicting colonies of red harvester ants (\emph{Pogonomyrmex
  barbatus}) results in non-overlapping dominance regions that resemble
Voronoi partitions.
Non-overlapping dominance regions akin to centroidal Voronoi partitions are
well documented in~\cite{GWB:74,MR-MH:80,QD-VF-MG:99} for the mouthbreeder
fish (\emph{Tilapia mossambica}).
Territoriality behavior and competition among prides of African lions
(\emph{Panthera leo}) are discussed in~\cite{AM-CP:09}.
Overall, numerous animal species achieve territory partitioning without a
central coordinating entity and without synchronized communication, but
rather relying upon asynchronous accidental interactions and stigmergy. To
the best of our knowledge about biological and engineering multiagent
systems, asynchronous territory partitioning has been barely studied,
see~\cite{MR-MH:80} for introductory ideas about animal behavior, and
mathematical models and analysis are lacking.

To finalize the literature review, here is a synopsis of the mathematical
ideas from diverse disciplines that we bring to bear on deployment and
partitioning problems.  First, we adopt the so-called \emph{gossip
  communication model}, in which only peer-to-peer asynchronous
communication links are required. This communication model is widely
studied in the wireless communication literature; example references
include~\cite{DK-AD-JG:03,SB-AG-BP-DS:06}. 
Moreover, we consider control systems on a non-Euclidean state space; the
interest for non-Euclidean spaces has a rich history in nonlinear control
theory, dynamical systems and robotics, including the early
work~\cite{RWB:72} and a recent application to multiagent
systems~\cite{AS-RS:07}. 
Finally, we adopt various tools from topology and from the study of
hyperspaces of sets~\cite{SBN:78}.

\subsection*{Statement of contributions}
This paper uncovers novel mathematical principles and tools of relevance in
coordination problems and multiagent systems.

We tackle partitioning and coverage control algorithms in innovative ways.
First, we design algorithms that require only gossip communication, i.e.,
asynchronous, pairwise, and possibly unreliable communication.  Gossip
communication is a simple, robust and effective protocol for noisy and
uncertain wireless environments.  Gossip communication may be implemented
in wandering robots with short-range unreliable communication (see the
illustrative motion coordination strategy in
Appendix~\ref{app:robotic-implementation}).  Second, we propose a paradigm
shift in coverage control and multicenter optimization.
Classically~\cite{JC-FB:08r,RG-JC:07,IIH-DMS:07,SM-JC-FB:04n,MP-EF-FB:08h,MS-DR-JJS:08,RV-DAWB-TG-TH-MF-JM:09},
the state space for the coverage algorithms are the agents' positions,
i.e., as a function of the agents' positions the environment is divided
into regions and regions are assigned to each agent. Instead, in our
approach, the agents' positions are no longer a concern: the state space is
the space of partitions of the environment and the algorithm dictates how
to update the regions.

Within the innovative context of gossip communication and partition-based
mechanisms, we devise a novel algorithm for multicenter and coverage
optimization.  Our gossip coverage algorithm is a peer-to-peer version of
Lloyd algorithm and aims to compute centroidal Voronoi partitions. Which
robot pair communicates at any given time is the outcome of either a
deterministic or a stochastic process.  We also propose a modified version that
restricts communication exchanges to adjacent regions and that has suitable
continuity properties. Simulations illustrate that our algorithms
successfully compute centroidal Voronoi partitions.

To formally establish the convergence properties of our proposed gossip
coverage algorithms, we perform a detailed mathematical analysis composed
of three steps.
First, we develop suitable versions of the Krasovskii-LaSalle invariance
principle for dynamical systems on metric spaces described by deterministic
and stochastic switches.  Convergence to a set of fixed points is achieved
under a uniform deterministic or stochastic persistency condition.
Second, we establish the continuity of various geometric maps, including
(1) the Voronoi diagram as a function of its generators, (2) the location
of the generalized centroid as a function of a set, (3) the widely-known
multicenter function studied in facility location, and (4) the gossip
coverage algorithms. These continuity properties are established with
respect to the symmetric distance metric in the space of partitions.
Third and final, we study the topology of the space of partitions.  With
respect to the symmetric difference metric we prove the compactness of a
relevant subset of partitions. Specifically, we focus on partitions whose
component regions are the union of a bounded number of convex sets.

In summary, relying upon our extensions of the invariance principle, the
compactness of a subset of the set of partitions, and the continuity of the
various relevant maps, we establish the convergence properties of the
proposed gossip algorithms.  In short, the algorithms converge to the set
of centroidal Voronoi partitions under mild technical assumptions and under
the assumption that the gossip communication exchanges satisfy either a
deterministic or a stochastic persistency condition.

\subsection*{Organization and notations}
The paper is structured as follows. In
Section~\ref{sec:review+preliminaries} we review multicenter optimization
and coverage control ideas. In Section~\ref{sec:OurAlgo} we state our
asynchronous territory partitioning problem, provide a solution via the
gossip coverage algorithm, state the convergence properties of the
algorithm and report some simulation results.  The following
Sections~\ref{sec:Krasovskii-LaSalle},~\ref{sec:PartitionsSpace}
and~\ref{sec:ContinuityMaps} develop the mathematical machinery required to
prove the convergence results.  Section~\ref{sec:Krasovskii-LaSalle}
contains the convergence theorems extending the Krasovskii-LaSalle
invariance principle.  Section~\ref{sec:PartitionsSpace} contains a
discussion about the compactness properties of the space of partitions.
Section~\ref{sec:ContinuityMaps} states the continuity properties of the
relevant maps and functions and contains the proof of the main convergence
results. Concluding remarks are given in Section~\ref{sec:Outro}.

We let $\realpositive$ and $\realnonnegative$ denote the set of positive
and non-negative real numbers, respectively, and $\integernonnegative$
denote the set of non-negative integer numbers.  Given $A\subset\real^d$,
we let $\interior(A)$, $\closure{A}$, $\partial A$ and $\diam(A)$ denote
its interior, its closure, its boundary and its diameter, respectively.
Given two sets $X$ and $Y$, a \emph{set-valued map} $\setmap{T}{X}{Y}$
associates to an element of $X$ a subset of $Y$.

\section{A review of multicenter optimization and distributed coverage
  control}
\label{sec:review+preliminaries}
In this section we review a variety of known results in geometric
optimization and in robotic coordination.  In
Subsection~\ref{subsec:multicenter-optimization} we review the notion of
environment partitions and we introduce the multicenter function as a way
to define optimal environment partitions and optimal robot or sensor
positions in the environment. In
Subsection~\ref{subsec:distribcoverage+delaunaynetwork} we review some
distributed control algorithms for agent motion coordination and
environment partitioning based on the classic Lloyd algorithm.

\subsection{Partitions, centroids and multicenter optimization}
\label{subsec:multicenter-optimization}
We let $Q$ denote an environment of interest to be apportioned. We assume
$Q$ is a compact convex subset of $\real^d$ with non-empty interior.
Partitions of $Q$ are defined as follows.
\begin{definition}[Partition]
  \label{def:partition}
  An \emph{$N$-partition of $Q$}, denoted by $v=(v_i)_{i=1}^N$, is an
  ordered collection of $N$ subsets of $Q$ with the following properties:
  \begin{enumerate}
  \item\label{def:Cover} $\union_{i\in\until{N}} v_i = Q$;
  \item\label{def:NullInters} $\interior(v_i) \intersect \interior(v_j)$ is
    empty for all $i,j\in\until{N}$ with $i\not=j$; and
  \item\label{def:Non-EmptyPart} each set $v_i$, $i\in\until{N}$, is
    closed and has non-empty interior.
  \end{enumerate}
  The set of $N$-partitions of $Q$ is denoted by $\part$.
\end{definition}

Let $p=(p_1,\ldots,p_N)\in Q^N$ denote the position of $N$ agents in the
environment $Q$.  Given a group of $N$ agents and an $N$-partition, each
agent is naturally in one-to-one correspondence with a component of the
partition; specifically we refer to $v_i$ as the \emph{dominance region} of
agent $i\in\until{N}$.

On $Q$, we define a \emph{density function} to be a bounded measurable
positive function $\map{\phi}{Q}{\realpositive}$ and a \emph{performance
  function} to be a locally Lipschitz, monotone increasing and convex
function $\map{f}{\realnonnegative}{\realnonnegative}$.  With these
notions, we define the \emph{multicenter function}
$\map{\Hgeneric}{\part\times Q^N}{\realnonnegative}$ by
\begin{equation}
  \label{eq:H(v,p)}
  \Hgeneric(v,p) = \sum_{i=1}^{N}\int_{v_i} f(\|p_i-q\|) \phi(q)dq.
\end{equation}
This function is well-defined because closed sets are measurable.  We aim
to minimize $\Hgeneric$ with respect to both the partition~$v$ and the
locations~$p$.

\begin{remarks}[Locational optimization]
  As discussed in the introduction and in the survey~\cite{QD-VF-MG:99},
  the multicenter function $\Hgeneric$ has numerous interpretations. Here
  we review two applications entailing robotic networks.  First, in a
  vehicle routing and service delivery example~\cite{MP-EF-FB:08h}, given
  vehicles at locations $p_i$, assume that $f(\|p_i-q\|)$ is the cost
  incurred by agent $i$ to travel to service an event taking place at point
  $q$.  Events take place inside $Q$ with likelihood $\phi$.  Accordingly,
  $\Hgeneric$ quantifies the expected wait-time between event arrivals and
  agents servicing them.

  Second, in an environmental monitoring application~\cite{SM-JC-FB:04n},
  assume the robots aim to detect acoustic signals that originate and
  propagate isotropically in the environment.  Because of noise and loss of
  resolution, the ability to detect a sound source originating at a
  point~$q$ from a sensor at position $p_i$ is proportional to the
  signal-to-noise ratio (which degrades with $\|q-p_i\|$). If the
  performance function $f$ equals minus the signal-to-noise ratio, then
  $\Hgeneric$ quantifies the expected signal-to-noise ratio and detection
  capacity for acoustic signals generated at random locations.  \oprocend
\end{remarks}

Among all possible ways of partitioning a subset of $\real^d$, one is worth
of special attention.  Define the \emph{set of partly coincident locations}
$S_N=\setdef{p\in Q^N}{p_i=p_j\text{ for some } i,j\in\until{N},\,
  i\not=j}$.  Given $p\in Q^N\setminus S_N$, the \emph{Voronoi partition of
  $Q$ generated by $p$}, denoted by $V(p)$, is the ordered collection of
the \emph{Voronoi regions} $\big(V_i(p)\big)_{i=1}^N$, defined by
\begin{equation}
  \label{eq:DefVor}
    V_i(p) = \setdef{q\in Q} {\|q-p_i\|\le\|q-p_j\| \text{ for all } j\neq i }.
\end{equation}
In other words, the Voronoi partition is a map $\map{V}{(Q^N\setminus
  S_N)}{\part}$.  The regions $V_i(p)$, $i\in\until{N}$, are convex and, if
$Q$ is a polytope, they are polytopes.  Now, given two distinct points
$q_1$ and $q_2$ in $\real^d$, define the \emph{$(q_1;q_2)$-bisector
  half-space} by
\begin{equation}
  \label{def:bisector-half-space}
  \halfs(q_1;q_2) = \setdef{q\in\real^d}{\|q-q_1\|\leq\|q-q_2\|}.
\end{equation}
In other words, the set $\halfs(q_1;q_2)$ is the closed half-space
containing $q_1$ whose boundary is the hyperplane bisecting the segment
from $q_1$ to $q_2$.  Note that bisector subspaces satisfy
$\halfs(q_1;q_2)\not=\halfs(q_2;q_1)$ and that Voronoi partition of $Q$
satisfies $V_i(p_1,\dots,p_N) = Q \intersection \big( \intersection_{j\neq
  i} \halfs(p_i;p_j) \big)$.

Each region equipped with a density function possesses a point with a
special relationship with the multicenter function.  Define the scalar
\emph{$1$-center function} $\Hone$ by
\begin{equation}
  \label{eq:H1def}
  \Hone(p;A)=\int_A f(\|p-q\|)\phi(q)dq,
\end{equation}
where $p$ is any point in $Q$ and $A$ is a compact subset of $Q$.  Under
the stated assumptions on the performance function $f$, the function
$p\mapsto\Hone(p;A)$ is strictly convex in $p$, for any set $A$ with
positive measure (we postpone the proof to
Lemma~\ref{lem:H1cont}). Therefore, the function $p\mapsto\Hone(p;A)$ has a
unique minimum in the compact and convex set $Q$.  We define the
\emph{generalized centroid} of a compact set $A\subset{Q}$ with positive
measure by
\begin{equation}
  \label{eq:def-Cd}
  \Cd(A) = \argmin \setdef{\Hone(p;A)}{p\in Q}.
\end{equation}
In what follows, it is convenient to drop the word ``generalized,'' and to
denote by $\Cd(v)=(\Cd(v_1),\dots,\Cd(v_N))\in Q^N$ the vector of regions
centroids corresponding to a partition $v\in\part$.

\begin{remark}[Quadratic and linear performance functions]
  If the performance function is $f(x)=x^2$, then the global minimum of
  $\Hone$ is the \emph{centroid} (also called the center of mass) of $A$,
  defined by
  \begin{equation*}
    \Cd(A)= \Big(\int_{A}\phi(q)dq\Big)^{-1} \int_{A}q\phi(q)dq.
  \end{equation*}
  If the performance function is $f(x)=x$, then the global minimum of
  $\Hone$ is the \emph{median} (also called the Fermat--Weber center) of
  $A$. See~\cite[Chapter 2]{FB-JC-SM:09} for more details.\oprocend
\end{remark}

\begin{figure}[thb]
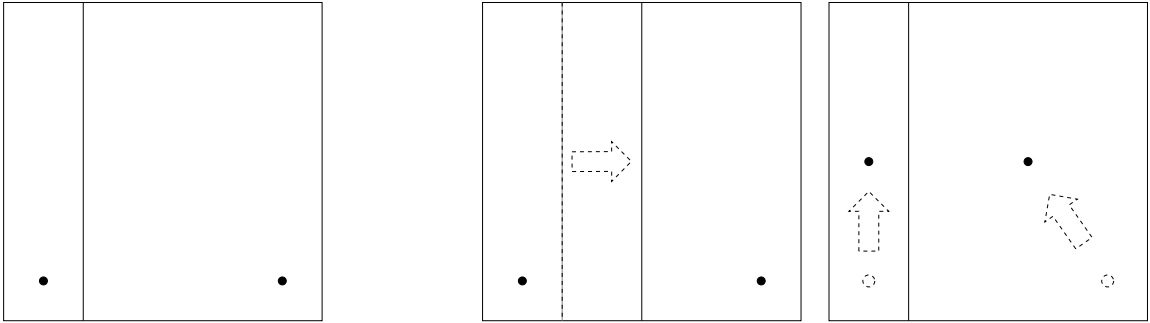
\centering%
  \begin{center}
    \resizebox{0.25\textwidth}{!}{\input{lloyd-illustration1.tex}}\qquad\qquad\qquad%
    \resizebox{0.25\textwidth}{!}{\input{lloyd-illustration3.tex}}\quad%
    \resizebox{0.25\textwidth}{!}{\input{lloyd-illustration2.tex}}%
  \end{center}
  \caption{Illustration of Proposition~\ref{prop:optimal-for-Hgeneric}.
    The left figure shows a sample $2$-partition $v$ and point set $p$ in a
    uniform square environment. The value of the cost function $\Hgeneric$
    at $(v,p)$ is diminished by either replacing $v$ with the Voronoi
    partition generated by $p$ (see central figure), or replacing $p$ with
    the centroids of $v$ (see right figure).  }
  \label{fig:optimize-p-or-v}
\end{figure}

Voronoi partitions and centroids have useful optimality properties stated
in the following proposition and illustrated in
Fig.~\ref{fig:optimize-p-or-v}.

\begin{proposition}[Properties of $\Hgeneric$]
  \label{prop:optimal-for-Hgeneric}
  For any partition $v\in \part$ and any point set $p\in Q^N\setminus S_N$,
  \begin{align}
    \Hgeneric(V(p),p) &\le \Hgeneric(v,p), \label{eq:Hgen-1}
    \\
    \Hgeneric(v,\Cd(v)) &\le \Hgeneric(v,p). \label{eq:Hgen-2}
  \end{align}
  Furthermore, inequality~\eqref{eq:Hgen-1} is strict if any entry of
  $V(p)$ differs from the corresponding entry of $v$ by a set with positive
  measure, and inequality~\eqref{eq:Hgen-2} is strict if $\Cd(v)$ differs
  from $p$.
\end{proposition}

The statements in Proposition~\ref{prop:optimal-for-Hgeneric} originate in
the early work by S.\ P.\ Lloyd~\cite{SPL:82}; modern treatments are given
in~\cite{QD-VF-MG:99} and~\cite[Propositions~2.14 and 2.15]{FB-JC-SM:09}.
Proposition~\ref{prop:optimal-for-Hgeneric} implies the following necessary
condition: if a pair $(v,p)$ with $p\not\in{S_N}$ minimizes $\Hgeneric$,
then $p=\Cd(v)$ and $v=V(p)$ up to a set of measure zero. Accordingly, the
partitions that minimize $\Hgeneric$ have the following property.
\begin{definition}\label{def:centroidal-Voronoi}
  The partition $v\in\part$ is \emph{centroidal Voronoi} if it has distinct
  centroids, that is, $\Cd(v_i)\neq\Cd(v_j)$ for all $j\neq i$, and
  $v=V(\Cd(v))$.
\end{definition}

\subsection{From Lloyd algorithm to distributed coverage control}
\label{subsec:distribcoverage+delaunaynetwork}
Here we consider a group of robotic agents with motion, communication and
computation capacities and we review a coverage control algorithm that
determines the motion of each robot in a group and the associated partition
in such a way as to minimize $\Hgeneric$.  In what follows, we restrict our
attention to $d=2$, that is, we assume $Q\subset\real^2$.

\begin{figure}[htb]
  \centering%
  \includegraphics[width=.4\linewidth]{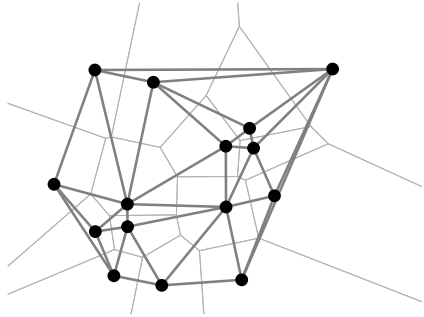}
  \caption{The Voronoi partition and the corresponding Delaunay graph}
  \label{fig:voronoi+delaunay}
\end{figure}

To explain in what sense our algorithms are distributed, we introduce a
useful graph.  The \emph{Delaunay
  graph}~\cite{MdB-MvK-MO-OS:00,FB-JC-SM:09} associated to the distinct
positions $p\in Q^N\setminus S_N$ is the undirected graph with node set
$\{p_i\}_{i=1}^N$ and with the following edges: $(p_i,p_j)$ is an edge if
and only if $V_i(p)\intersect V_j(p)$ is non-empty.  In other words, two
agents are neighbors if and only if their Voronoi regions intersect, see
Fig.~\ref{fig:voronoi+delaunay}.

The \emph{coverage algorithm} we consider is a distributed version of the
classic Lloyd algorithm~\cite{QD-VF-MG:99} based on ``centering and
partitioning'' for the computation of centroidal Voronoi partitions. The
algorithm is distributed in the sense that each robot determines its region
of responsibility and its motion plan based upon communication with only
some neighbors. Specifically, communications among the robots takes place
along the edges of the Delaunay graph.  The distributed coverage algorithm
is described as follows.  At each discrete time instant
$t\in\integernonnegative$, each agent $i$ performs the following tasks: (i)
it transmits its position and receives the positions of its neighbors in
the Delaunay graph; (ii) it computes its Voronoi region with the
information received; (iii) it moves to the centroid of its Voronoi region.
In mathematical terms, for $t\in\integernonnegative$,
\begin{equation}
  \label{eq:coverage-control}
  p(t+1)= \Cd(V(p(t))).
\end{equation}
A variation of the function $\Hgeneric$ is useful to analyze this
algorithm. We define the positions-based multicenter function
$\map{\Hvoronoi}{Q^N\setminus S_N}{\realnonnegative}$ by
\begin{align}   \label{eq:def-Hvoronoi}
  \Hvoronoi(p)&=\Hgeneric(V(p),p)
  =\sum_{i=1}^{N}\int_{V_i(p)} f(\|q-p_i\|) \phi(q)d q.
\end{align}
Because of the compactness of $Q$, a continuity property, and the
monotonicity properties in Proposition~\ref{prop:optimal-for-Hgeneric}, one
can show~\cite[Theorem~5.5]{FB-JC-SM:09} that $\Hvoronoi$ is monotonically
non-increasing along the solutions of~\eqref{eq:coverage-control} and that
all solutions of~\eqref{eq:coverage-control} converge asymptotically to the
set of configurations that generate centroidal Voronoi
partitions. Additional considerations about convergence are given
in~\cite{QD-ME-JL:06}.


\section{Gossip coverage control as a dynamical system on the space of
  partitions}
\label{sec:OurAlgo}
In this section we present the problem of interest, our novel gossip
coverage algorithm and its convergence properties in
Subsections~\ref{subsec:problem-statement},
\ref{subsec:gossip-coverage-algo} and~\ref{sec:AlgoConverg}, respectively.
In order to reduce the communication requirements of our algorithm, we
propose an adjacency-based and continuous algorithm in
Subsection~\ref{sec:SmoothedMap}.  Finally, we report some simulation
results in Subsection~\ref{subsec:simulation-results}.

\subsection{Problem statement}
\label{subsec:problem-statement}

The distributed coverage law, based upon the Lloyd algorithm and described
in the previous section, has some important limitations: it is applicable
only to robotic networks with \emph{synchronized} and \emph{reliable}
communication along \emph{all edges of the Delaunay graph} (computed as a
function of the robots' positions).  In other words, the
law~\eqref{eq:coverage-control} requires that there exists a predetermined
common communication schedule for all robots and, at each communication
round, each robot must simultaneously and reliably communicate its
position.  Note that the Delaunay graph, interpreted as a communication
graph, has the following drawbacks: for worst-case robots' positions, a
robot might have $N-1$ neighbors in the Delaunay graph and/or might have a
neighbor that is arbitrarily far inside the environment. Therefore, each
robot must be capable to communicate potentially to all other robots and/or
to robots at large distances.

Given this broad range of undesirable limitations, the aim of this paper is
to reduce the communication requirements of distributed coverage
algorithms, in terms of reliability, synchronization and topology. Here are
some relevant questions that constitute our informal problem statement:
\begin{quote}
  Is it possible to optimize robots positions and environment partition
  with asynchronous, unreliable, and delayed communication?  Specifically,
  what if the communication model is that of \emph{gossiping agents}, that
  is, a model in which only a pair of robots can communicate at any time?
  Since Voronoi partitions generated by gossiping and moving agents cannot
  be computed by gossiping agents, how do we update the environment
  partition?
\end{quote}
To answer these questions, the next subsections propose an innovative
partition-based gossip approach, in which the robots' positions essentially
play no role and where instead dominance regions are iteratively updated.
Designing coverage algorithms as dynamical systems on the space of
partitions has the key advantage that one is not restricted to working only
with Voronoi or anyway position-dependent partitions.

\begin{example}[The Lloyd algorithm in the partition-based approach] 
  The distributed coverage algorithm~\eqref{eq:coverage-control} updates
  the robots' positions so as to incrementally minimize the function
  $\Hvoronoi$, while the environment partition is a function of the robots'
  positions.  In this paper we take a dual approach: we consider an
  algorithm that evolves partitions.  From this partition-based viewpoint,
  the coverage algorithm is an iterated map on $\part$ and
  equation~\eqref{eq:coverage-control} is rewritten as
  $v(t+1)=V(\Cd(v(t)))$.  \oprocend
\end{example}

\subsection{The  gossip coverage algorithm}
\label{subsec:gossip-coverage-algo}
In this subsection we present a novel partition-based coverage algorithm in
which, at each communication round, only two regions communicate.
Recall the notion of bisector half-space from
equation~\eqref{def:bisector-half-space}.

\renewcommand\footnoterule{\hrule width \textwidth height .4pt}

\medskip\footnoterule\vspace{.75\smallskipamount}
\noindent\hfill\textbf{Gossip Coverage Algorithm}\hfill\vspace{.75\smallskipamount}
\footnoterule\vspace{.75\smallskipamount}
\noindent For all $t\in\integernonnegative$, each agent $i\in\until{N}$
maintains in memory a dominance region $v_i(t)$.  The collection
$(v_1(0),\dots,v_N(0))$ is an arbitrary polygonal $N$-partition of $Q$.  At
each $t\in\integernonnegative$ a pair of communicating regions, say
$v_i(t)$ and $v_j(t)$, is selected by a deterministic or stochastic process
to be determined.  Every agent $k\not\in\{i,j\}$ sets $v_k(t+1):=v_k(t)$.
Agents $i$ and $j$ perform the following tasks:
\begin{algorithmic}[1]
  \STATE agent $i$ transmits to agent $j$ its dominance region $v_i(t)$ and
  vice-versa%
  \STATE both agents compute the centroids $\Cd(v_i(t))$ and
  $\Cd(v_j(t))$ %
  \IF{$\Cd(v_i(t))= \Cd(v_j(t))$} %
  \STATE $v_i(t+1) := v_i(t)$ and $v_j(t+1) := v_j(t)$ \ELSE %
  \STATE $v_i(t+1) := \big(v_i(t)\union v_j(t)\big) \intersect
  \halfs\big(\!\Cd(v_i(t));\Cd(v_j(t))\big)$ \\
  $v_j(t+1) := \big(v_i(t)\union v_j(t)\big) \intersect
  \halfs\big(\!\Cd(v_j(t));\Cd(v_i(t))\big)$
  \ENDIF
\end{algorithmic}
\vspace{.5\smallskipamount}\footnoterule\smallskip

\begin{figure}[thb]
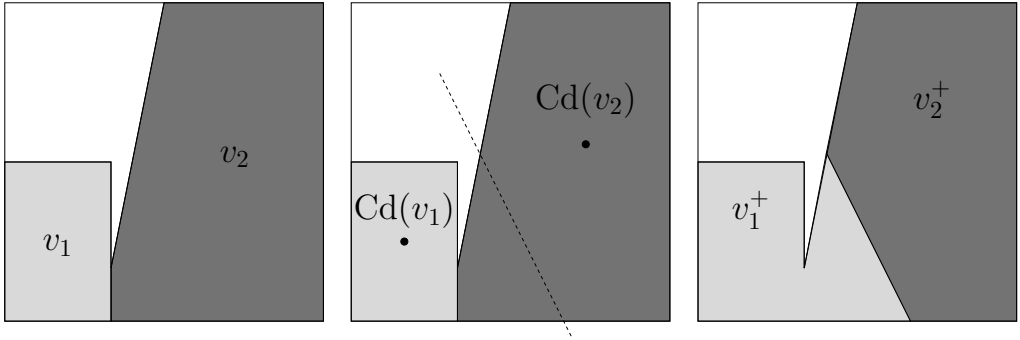
\centering%
  \begin{center}
    \resizebox{0.25\textwidth}{!}{\input{gc1.tex}}\quad%
    \resizebox{0.25\textwidth}{!}{\input{gc2.tex}}\quad%
    \resizebox{0.25\textwidth}{!}{\input{gc3.tex}}%
  \end{center}
  \caption{The gossip coverage algorithm. The left and right figure contain
    the initial partition and the partition after one application of the
    gossip coverage algorithm.  In the middle figure we show the two
    centroids and (with a dashed line) the segment determining the bisector
    half-space.}
  \label{fig:ExamplePairwise-simple}
\end{figure}

In other words, when two agents with distinct centroids communicate, their
dominance regions evolve as follows: the union of the two dominance regions
is divided into two new dominance regions by the hyperplane bisecting the
segment between the two centroids; see
Fig.~\ref{fig:ExamplePairwise-simple}.  As a consequence, if the
centroids $\Cd(v_i(t))$, $\Cd(v_j(t))$ are distinct, then
$\{v_i(t+1),v_j(t+1)\}$ is the Voronoi partition of the set $v_i(t)\union
v_j(t)$ generated by the centroids $\Cd(v_i(t))$ and $\Cd(v_j(t))$.

We claim that the algorithm is well-posed in the sense that the sequence of
collections $\{v(t)\}_{t\in\integersnonnegative}$ generated by the
algorithm is an $N$-partition at all times $t$. Indeed, it is immediate to
see that the first two properties in Definition~\ref{def:partition} are
satisfied at all time if they are satisfied at initial time.  Finally, at
all times $t$, each component of $v(t)$ is clearly closed and has non-empty
interior because of the following geometric fact: there exists no
half-plane containing the interior of a region and not containing the
centroid of the same region.

Now, for any $i,j\in\until{N}$ with $i\not=j$, define the map
$\map{T_{ij}}{\part}{\part}$ by
\begin{equation*}
  T_{ij}(v) =
  \begin{cases}
    v,   &\text{if } \Cd(v_i)=\Cd(v_j),\\
    (v_1,\dots,\widehat{v}_i,\dots,\widehat{v}_j, \dots,v_N),
    \quad &\text{otherwise},
  \end{cases}
\end{equation*}
where
\begin{equation}
  \label{eq:lastminute}
  \begin{split}
  \widehat{v}_i &= \big(v_i\union v_j\big) \intersect
  \halfs\big(\!\Cd(v_i);\Cd(v_j)\big),\\
  \widehat{v}_j &= \big(v_i\union v_j\big) \intersect
  \halfs\big(\!\Cd(v_j);\Cd(v_i)\big).
  \end{split}
\end{equation}
The dynamical system on the space of partitions is therefore described by,
for $t\in\integernonnegative$,
\begin{equation}
  \label{eq:OurAlgo-ij}
  v(t+1)=T_{ij}(v(t)),
\end{equation}
together with a rule describing what pair of regions $(i,j)$ is selected at
each time.  We also define the set-valued map $\setmap{T}{\part}{\part}$ by
$T(v)=\setdef{T_{ij}(v)}{i,j\in\until{N},\; i\neq j}$. The map $T$
describes one iteration of the gossip coverage algorithm; an evolution of
the gossip coverage algorithm is one of the solutions to the
non-deterministic set-valued dynamical system $v(t+1)\in T(v(t))$.

\begin{remark}[Gossip disk-covering control]
  \label{rem:disk-covering}
  We believe that our gossip and partition-based algorithmic approach is
  applicable to a broad range of coverage problems. For example, the
  worst-case multicenter function~\cite{JC-FB:08r} is defined by
  $\HHdc(p,v) = \max_{i\in\until{N}} \max_{q\in{v_i}} \|q-p_i\|$.
  Maximizing $\HHdc$ is equivalent to covering $Q$ with $N$ disks of
  smallest radius centered at $p$.  As in
  Proposition~\ref{prop:optimal-for-Hgeneric}, for any $v\in \part$ and
  $p\in Q^N\setminus S_N$, one can prove $\HHdc(V(p),p) \le \HHdc(v,p)$ and
  $\HHdc(v,\CC(v)) \le \HHdc(v,p)$, where $\CC(v)$ is the array of
  circumcenters of the components of $v$. Hence, a gossip coverage
  algorithm for $\HHdc$ is designed by replacing centroid with circumcenter
  operations in~\eqref{eq:lastminute}. We leave this and further extensions
  to future works. \oprocend
\end{remark}

\subsection{Analysis results for the gossip coverage algorithm}
\label{sec:AlgoConverg}
In this subsection we state the main analysis and convergence results for
the gossip coverage algorithm.

We begin by studying the fixed points of $T$ and by
  introducing an appropriate cost function with monotonicity properties
  along $T$.  Regarding the algorithm's fixed points, we extend
  Definition~\ref{def:centroidal-Voronoi} as follows.  A partition $v$ is
  \emph{mixed centroidal Voronoi} if, for all pairs $(v_i,v_j)$ with
  $i\not=j$, either $\Cd(v_i)=\Cd(v_j)$ or $(v_i,v_j)$ is a centroidal
  Voronoi partition of $v_i\union v_j$, that is, $v_i = (v_i\union v_j)
  \intersect \halfs\big(\!\Cd(v_i);\Cd(v_j)\big)$.
  \begin{lemma}[Fixed points of $T$ and centroidal Voronoi partitions]
    \label{lem:fixed-points-cvp}
    For $i,j\in\until{N}$, $j\neq i$, denote the set of fixed points of
    $\map{T_{ij}}{\part}{\part}$ by
    $F_{ij}=\setdef{v\in\part}{T_{ij}(v)=v}$. The following statements
    hold:
    \begin{enumerate}
    \item $\intersect_{j\neq i} F_{ij}$ equals the set of mixed centroidal
      Voronoi partitions; and
    \item\label{item:from-mixed-to-cvt} if $v$ is a mixed centroidal
      Voronoi partition satisfying $\Cd(v_i)\neq\Cd(v_j)$ for $j\neq i$,
      then $v$ is centroidal Voronoi.
    \end{enumerate}
  \end{lemma}%
  \noindent Next, we define the partition-based multicenter function
  $\map{\Hcentroid}{\part}{\realnonnegative}$ by
  \begin{align}
    \Hcentroid(v)&= \Hgeneric(v,\Cd(v)) =
    \sum_{i=1}^{N}\int_{v_i} f(\|q-\Cd(v_i)\|) \phi(q) dq.
    \label{eq:H-centroid}
  \end{align}
  \begin{lemma}[Monotonicity of $\Hcentroid$ along $T$]
    \label{lem:HcentroidPairwiseDecr}
    For $i,j\in\until{N}$, $i\not=j$,
    \begin{align*}
      \Hcentroid(T_{ij}(v))&\le \Hcentroid(v), \enspace &&\text{for all }
      v\in\part, \text{ and}\\ 
      \Hcentroid(T_{ij}(v))& < \Hcentroid(v), &&\text{iff $T_{ij}(v)$ and
        $v$ differ by a set of measure zero}.
    \end{align*}
  \end{lemma}%
  \noindent{The} proofs of Lemmas~\ref{lem:fixed-points-cvp}
  and~\ref{lem:HcentroidPairwiseDecr} consist of elementary manipulations
  and are omitted in the interest of brevity.  In short, we have
  established that the function $\Hcentroid$ monotonically decreases along
  each $T_{ij}$ when away from fixed points, and that centroidal Voronoi
  partitions are the fixed points of all $T_{ij}$ provided centroids are
  distinct.

We now prepare to state the main convergence result for $T$.  We need to
introduce some useful properties for sequences of partitions and for
switching signals.

\begin{definition}[Non-degeneracy]
  A sequence of $N$-partitions $\seqdef{v(t)}{t\in\integernonnegative}$ is
  \begin{enumerate}
  \item \emph{(uniformly) distinct centroidal} if there exists $\epsilon>0$
    such that, for all $t\in\integernonnegative$ and $i,j\in\until{N}$,
    $i\not=j$, one has $\|\Cd(v_i(t))-\Cd(v_j(t))\|\geq\eps$;

  \item \emph{(uniformly componentwise) non-vanishing} if there exists
    $\epsilon>0$ such that, for all $t\in\integernonnegative$ and
    $i\in\until{N}$, the Lebesgue measure of $v_i(t)$ is greater than
    $\epsilon$; and

  \item \emph{(uniformly componentwise) finitely convex} if there exists
    $\ell\in\natural$ such that, for all $t\in\integernonnegative$ and
    $i\in\until{N}$, the set $v_i(t)$ is the union of at most $\ell$ convex
    sets.
  \end{enumerate}
  Moreover, the sequence $v$ is said to be \emph{non-degenerate} if it is
  distinct centroidal, non-vanishing and finitely convex.
\end{definition}

For example, a sequence of partitions is finitely-convex if each component
of each partition in the sequence is the union of a uniformly bounded
number of polygons with a uniformly bounded number of vertices.

\begin{definition}[Uniform and random persistency]
  \label{def:persistency}
  Let $X$ be a finite set.
  \begin{enumerate}
  \item A map $\map{\sigma}{\integernonnegative}{X}$ is \emph{uniformly
      persistent} if there exists a duration $\Delta\in\natural$ such that,
    for each $x\in{X}$, there exists an increasing sequence of times
    $\seqdef{t_k}{k\in\integernonnegative}\subset\integernonnegative$
    satisfying $t_{k+1}-t_k\leq\Delta$ and $\sigma(t_k)=x$ for all
    $k\in\integernonnegative$.
  \item A stochastic process $\map{\sigma}{\integernonnegative}{X}$ is
    \emph{randomly persistent} if there exists a probability $p\in{]0,1[}$
    and a duration $\Delta\in\natural$ such that, for each $x\in{X}$ and
    for each $t\in\integernonnegative$, there exists $k\in\until{\Delta}$
    satisfying
    \begin{equation*}
      \Prob\big[\sigma(t+k)= x \,|\,
      \sigma(t),\dots,\sigma(1)\big]  \geq p.
    \end{equation*}
  \end{enumerate}
\end{definition}

We are now ready to state the main deterministic and stochastic convergence
results for the gossip coverage algorithm. It is convenient to postpone to
Section~\ref{sec:ConvergProofs} the theorem proof and the definition of
convergence in the space of partitions.

\begin{theorem}[Convergence under persistent gossip]
  \label{th:T-persistent-gossip}
  Consider the gossip coverage algorithm $T$ and the evolutions
  $\seqdef{v(t)}{t\in\integernonnegative} \subset \part$ defined by
  \begin{equation*}
    v(t+1) = T_{\sigma(t)}(v(t)), \quad \text{for } t\in\integernonnegative,
  \end{equation*}
  where $\map{\sigma}{\integernonnegative}{\myset}$ is either a deterministic
  map or a stochastic process.  Then the following statements hold:
  \begin{enumerate}
  \item if $\sigma$ is a uniformly persistent map, then each non-degenerate
    evolution $v$ converges to the set of centroidal Voronoi partitions;
    and

  \item if $\sigma$ is a randomly persistent stochastic process, then each
    evolution $v$, conditioned upon being non-degenerate, converges almost
    surely to the set of centroidal Voronoi partitions.
  \end{enumerate}
\end{theorem}

The statements of Theorem~\ref{th:T-persistent-gossip} rely upon the
assumption of non-degenerate evolutions. It is our conjecture that,
starting from generic polygonal partition, this assumption is typically
satisfied. We will discuss some numerical evidence to this effect in the
next section.

Lemma~\ref{lem:HcentroidPairwiseDecr} indicates how the function
$\Hcentroid$ plays the role of a Lyapunov function for the dynamical system
defined by $T$.  To provide a complete Lyapunov convergence proof of
Theorem~\ref{th:T-persistent-gossip}, we set out to establish three sets of
relevant results. First, we need to establish extensions of the
Krasovskii-LaSalle invariance principle for set-valued dynamical systems
over compact metric spaces.  Second, we need to establish the compactness
properties of the space of non-degenerate partitions. Third, we need to
establish the continuity of the relevant geometric maps.  These three
topics are the subjects of Section~\ref{sec:Krasovskii-LaSalle},
\ref{sec:PartitionsSpace} and~\ref{sec:ContinuityMaps}, respectively.

\subsection{Designing an adjacency-based and continuous algorithm}
\label{sec:SmoothedMap}%
The gossip coverage map $T$ has the undesirable feature that it entails
communication exchanges between any two regions.  This communication
requirement might be too onerous in some multiagent applications. Ideally
we would like to require communications only between \emph{adjacent}
regions, that is, regions whose boundaries touch, or between ``nearby''
regions.  We believe such communication requirements may be easily achieved
in robotic networks and we detail a sample implementation for robots with
range-dependent Poisson communication in
Appendix~\ref{app:robotic-implementation}.
Additionally, we aim to design a continuous gossip coverage map. We require
the modified map to be continuous for technical reasons: the invariance
principles we adopt for the convergence analysis require continuity of the
dynamical system.

Motivated by this discussion, we modify the map $T$ to rely upon only
adjacency-based communication and to be continuous.  First, we introduce a
pseudodistance notion between sets. Given closed $A\subset{Q}$ and
$B\subset{Q}$ with non-empty interior, define
\begin{equation*}
  \pdist(A,B) = \inf\setdef{\|a-b\|}{(a,b)\in \interior(A)\times\interior(B)}.
\end{equation*}
Second, we select a positive constant $\delta\ll\diam(Q)$ and denote by
$\map{\Tdelta}{\part}{\part}$ the \emph{modified gossip coverage map} to be
defined in what follows.  For any $i,j\in\until{N}$, $i\neq{j}$, we give
the following partial definition:
\begin{equation}\label{def:partial}
  \Tdelta_{ij}(v) =
  \begin{cases}
    v, & \text{if } \big( \|\Cd(v_i)-\Cd(v_j)\|=0 \big) \text{
      or\phantom{and}\hspace{-1em} }
    \big(\pdist(v_i,v_j)\geq\delta\big),\\
    T(v), & \text{if } \big( \|\Cd(v_i)-\Cd(v_j)\|\geq\delta \big) \text{
      and\phantom{or}\hspace{-1em} } \big(\pdist(v_i,v_j)=0\big).
  \end{cases}
\end{equation}
Therefore, if either $\Cd(v_i)$ and $\Cd(v_j)$ coincide or the
pseudodistance between $v_i$ and $v_j$ is larger than $\delta$, then
$\Tdelta_{ij}(v)=v$, that is, the map $\Tdelta_{ij}$ leaves the partition
unchanged.  Additionally, if the pseudodistance between the regions $v_i$
and $v_j$ is zero ($v_i$ and $v_j$ are adjacent) and the distance between
$\Cd(v_i)$ and $\Cd(v_j)$ is larger than $\delta$, then
$\Tdelta_{ij}(v)=T_{ij}(v)$.

\begin{figure*}[th]
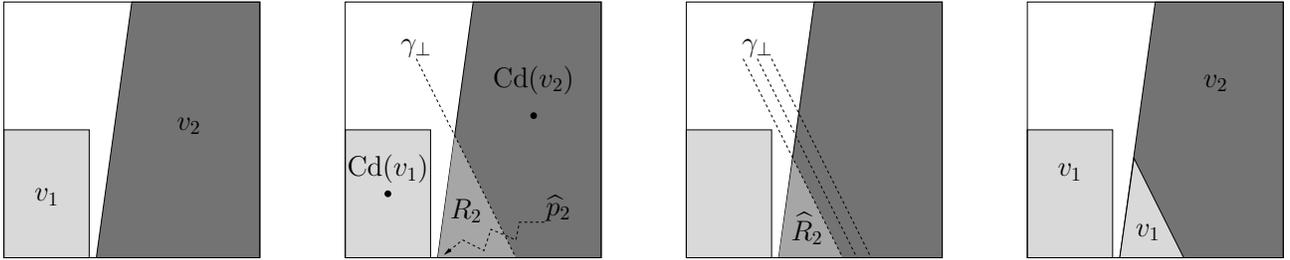
\centering%
  \begin{center}
    \resizebox{!}{0.15\textheight}{\input{gcs1.tex}}\hfill%
    \resizebox{!}{0.15\textheight}{\input{gcs2.tex}}\hfill%
    \resizebox{!}{0.15\textheight}{\input{gcs3.tex}}\hfill%
    \resizebox{!}{0.15\textheight}{\input{gcs4.tex}}
  \end{center}
  \caption{Modified gossip between close but not adjacent regions
    ($0<\beta_{12}(v)<1$).  The bisecting line~$\gamma_\perp$ borders the
    set $R_2=v_2 \intersect \halfs(\Cd(v_1),\Cd(v_2))$ that, in the map
    $T$, is assigned to $v_1$ (see Fig.~\ref{fig:ExamplePairwise-simple}).
    According to~$\Tdelta$ instead, only the subset $\widehat{R}_2
    \subsetneq R_2$ is assigned to $v_1$. Loosely speaking, the ``width''
    of $\widehat{R}_2$ equals $\beta_{12}(v)$ times the ``width'' of $R_2$,
    where ``width'' of a set is the maximum distance from a point in the
    set to $\gamma_\perp$.  }
\label{fig:TdeltaDemo}
\end{figure*}

Next, we consider partitions that do not satisfy either of the two
conditions in definition~\eqref{def:partial}.  We define the \emph{unit
  saturation function} $\map{\sat}{\realnonnegative}{[0,1]}$ by $\sat(x) =
x$ if $x\in[0,1]$, and $\sat(x) = 1$ if $x>1$ and the scaling function
$\map{\beta_{ij}}{\part}{[0,1]}$ by
\begin{equation*}
  \beta_{ij}(v) =\sat\!\big( \|\Cd(v_i)-\Cd(v_j)\|/\delta\big)
  \Big( 1 - \sat\!\big(\pdist(v_i,v_j)/\delta\big)  \Big).
\end{equation*}
The first condition and the second condition in~\eqref{def:partial}
correspond precisely to $\beta_{ij}(v)=0$ and $\beta_{ij}(v)=1$,
respectively.  For partitions $v$ satisfying $0<\beta_{ij}(v)<1$, we aim to
define $\Tdelta$ so as to continuously interpolate between the identity map
and the map $T$; see Fig.~\ref{fig:TdeltaDemo} for an illustration.
Let $R_i= v_i \intersect \halfs(\Cd(v_j),\Cd(v_i))$ and $R_j= v_j
\intersect \halfs(\Cd(v_i),\Cd(v_j))$. Define the line
$\gamma_\perp=\partial\halfs(\Cd(v_j),\Cd(v_i))$ and
\begin{equation}\label{eq:pi-pj}
  \begin{split}
    \widehat{p}_i & = \text{a point in } \closure{\interior{(R_i)}} \text{
      that is maximally distant from }\gamma_{\perp},\\
    \widehat{p}_j & = \text{a point in } \closure{\interior{(R_j)}} \text{
      that is maximally distant from }\gamma_{\perp}.
  \end{split}
\end{equation}
Next, note that for each $q\in{R_i\cup{R_j}}$ there exists a unique line,
say $\gamma_q$, that is parallel to $\gamma_\perp$ and passes through $q$.
Based on this notion, we define
\begin{align*}
  \widehat{R}_i &= \setdef{q\in R_i} { \dist(\widehat{p}_i,\gamma_q)\leq
    \beta_{ij}(v) \dist(\widehat{p}_i,\gamma_{\perp})
    \,\,\,\text{or}\,\,\, \dist(q,
    \gamma_{\perp})\geq \dist(\widehat{p}_i,\gamma_{\perp}) },\\
  \widehat{R}_j &= \setdef{q\in R_j} { \dist(\widehat{p}_j,\gamma_q)\leq
    \beta_{ij}(v) \dist(\widehat{p}_j,\gamma_{\perp})\,\,\,\text{or}\,\,\,
    \dist(q, \gamma_{\perp})\geq \dist(\widehat{p}_j,\gamma_{\perp}) }.
\end{align*}
We can now complete the partial definition~\eqref{def:partial}. For all $v$
with $0<\beta_{ij}(v)<1$, that is, for all partitions not already dealt with
in definition~\eqref{def:partial}, we define
\begin{equation*}
  \Tdelta_{ij}(v) =
  (v_1,\dots,\underbrace{\closure{\big(v_i\setminus\widehat{R}_i\big) \cup
      \widehat{R}_j }}_{i\text{th entry}} ,
  \dots, \underbrace{\closure{\big(v_j\setminus\widehat{R}_j\big) \cup
      \widehat{R}_i}}_{j\text{th entry}}, \dots,v_N).
\end{equation*}
As discussed for $T$, one can prove that the map
$\setmap{\Tdelta}{\part}{\part}$ defined by $ \Tdelta(v) =
\setdef{\Tdelta_{ij}(v)} {i,j\in\until{N},\,i\neq j}$, is well-posed and
has the following properties.

\begin{theorem}[Convergence of modified gossip map]
  \label{th:Tdelta-persistent-gossip}
  Consider the modified gossip coverage algorithm $\Tdelta$ and the
  evolutions $\seqdef{v(t)}{t\in\integernonnegative}\subset\part$ defined
  by
  \begin{equation*}
    v(t+1) = \Tdelta_{\sigma(t)}(v(t)), \quad \text{for } t\in\integernonnegative,
  \end{equation*}
  where $\map{\sigma}{\integernonnegative}{\myset}$ is either a deterministic
  map or a stochastic process.  Then the following statements hold:
  \begin{enumerate}
  \item if $\sigma$ is a uniformly persistent map, then each non-vanishing and
    finitely-convex evolution $v$ converges to the set of mixed centroidal
    Voronoi partitions; and

  \item if $\sigma$ is a randomly persistent stochastic process, then each
    evolution $v$, conditioned upon being non-vanishing and finitely
    convex, converges almost surely to the set of mixed centroidal Voronoi
    partitions.
  \end{enumerate}
\end{theorem}

\subsection{Simulation results and implementation remarks}
\label{subsec:simulation-results}
\begin{figure*}[th]\centering%
\includegraphics[width=.199\linewidth]{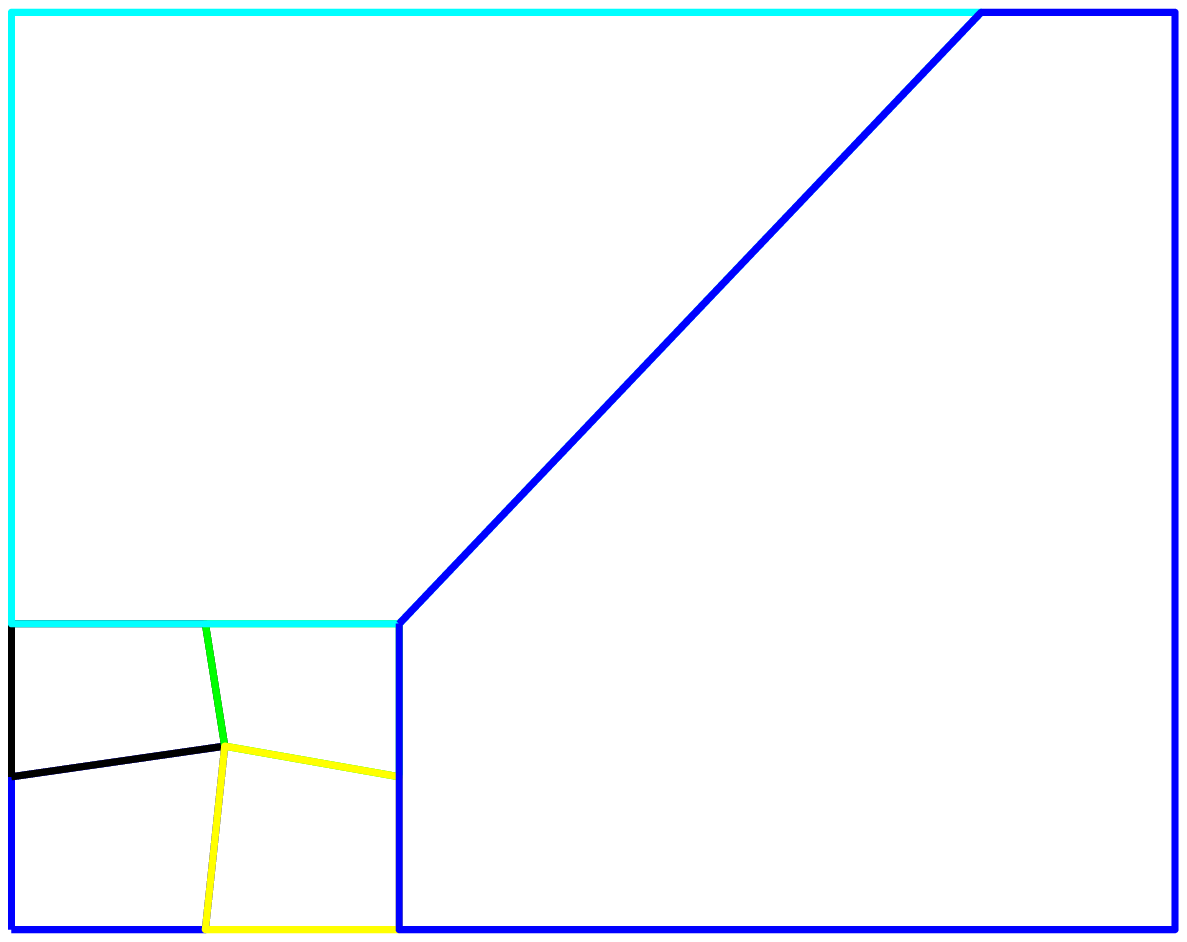}%
\includegraphics[width=.199\linewidth]{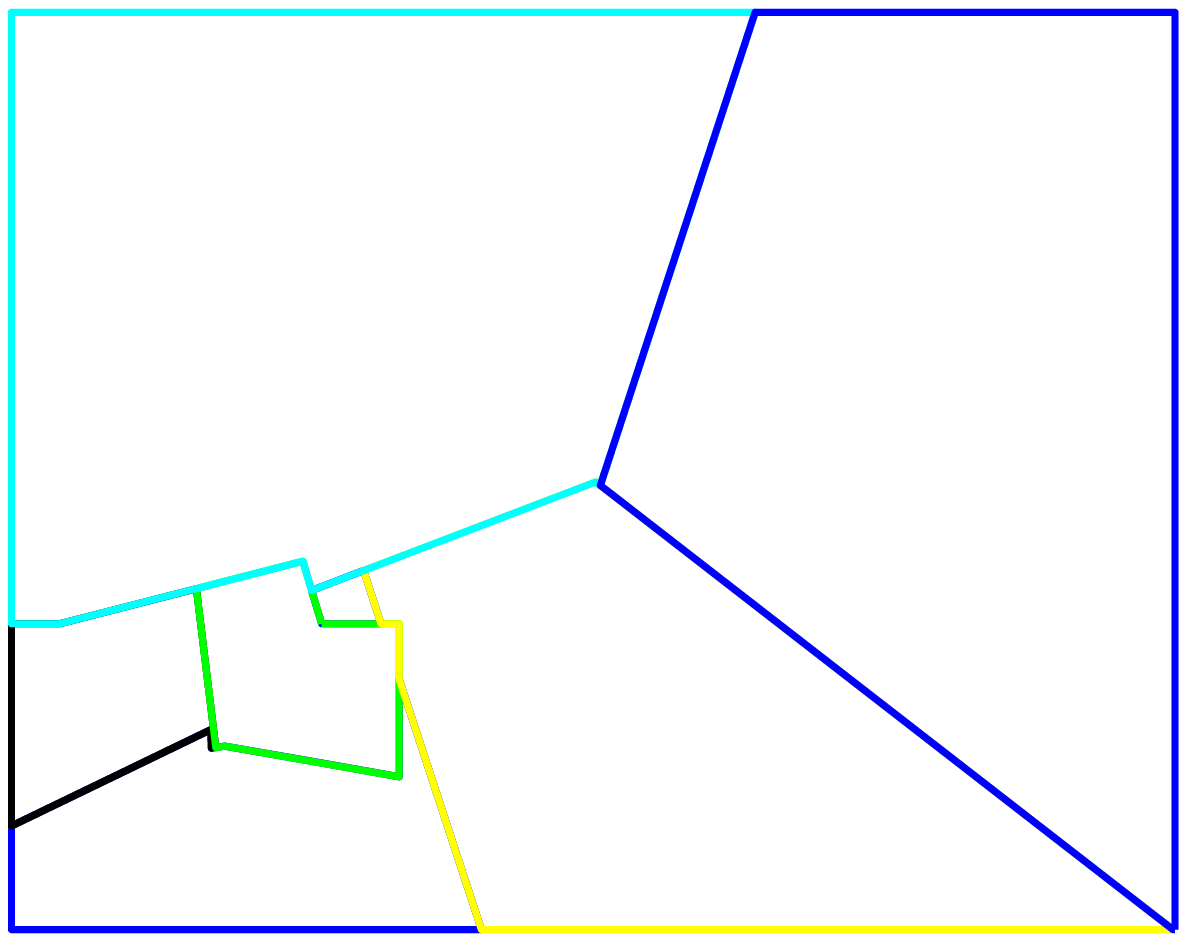}%
\includegraphics[width=.199\linewidth]{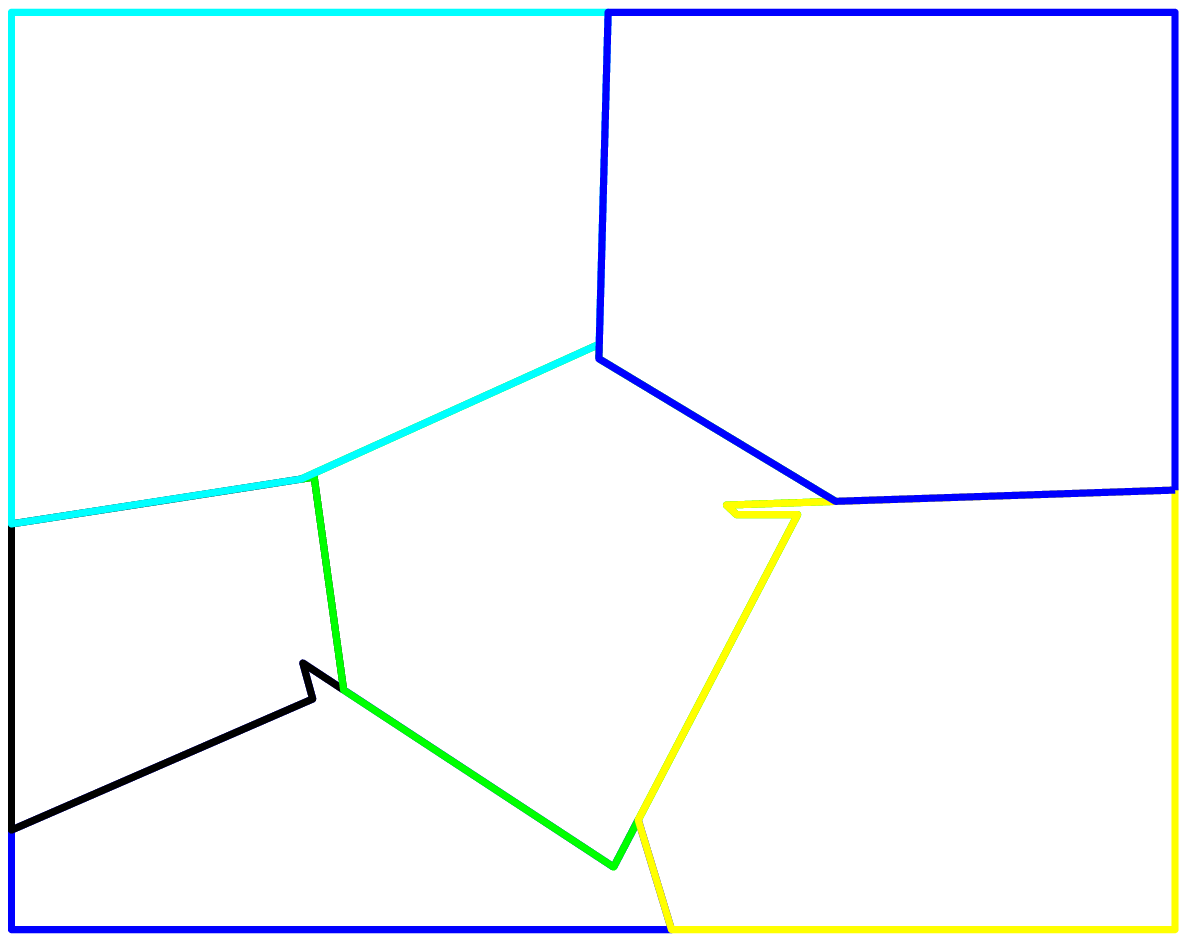}%
\includegraphics[width=.199\linewidth]{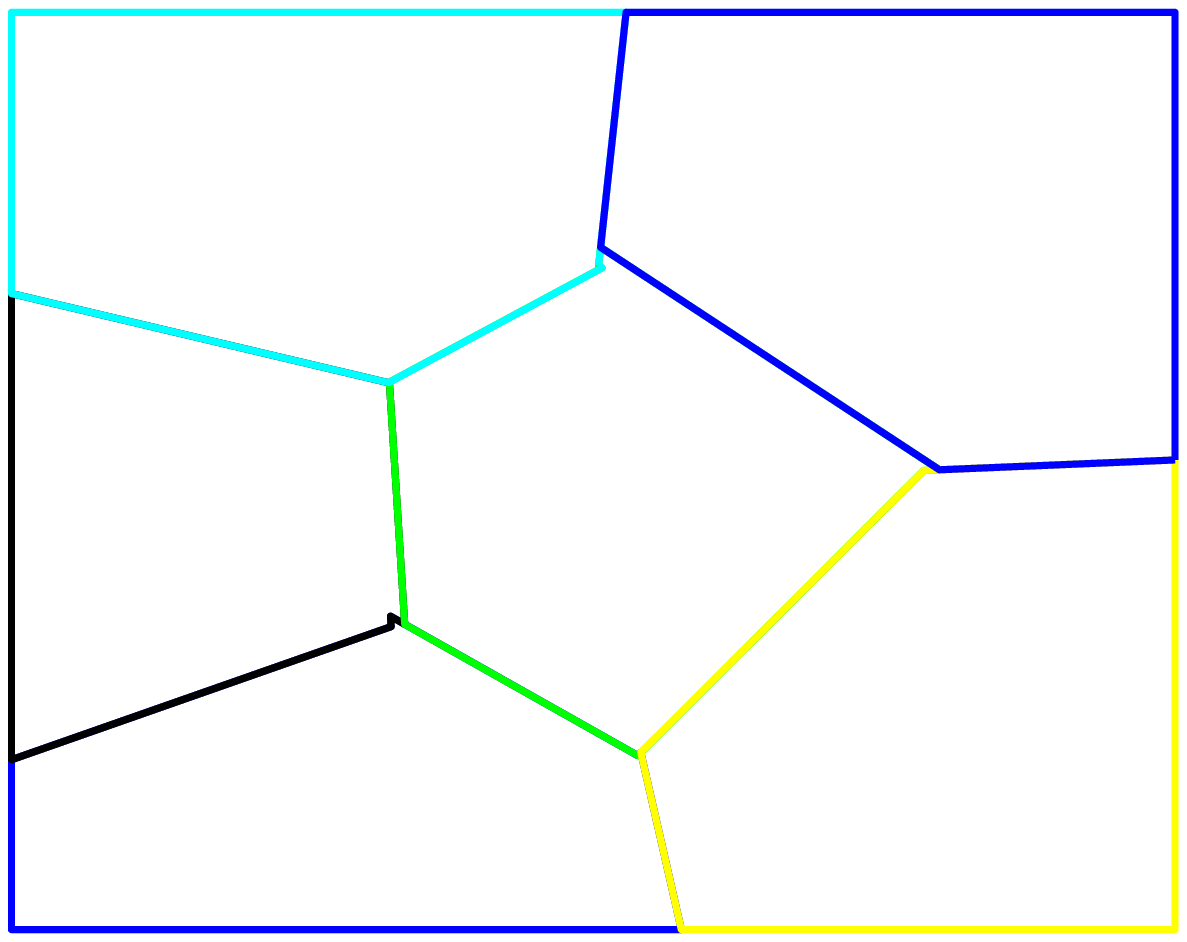}%
\includegraphics[width=.199\linewidth]{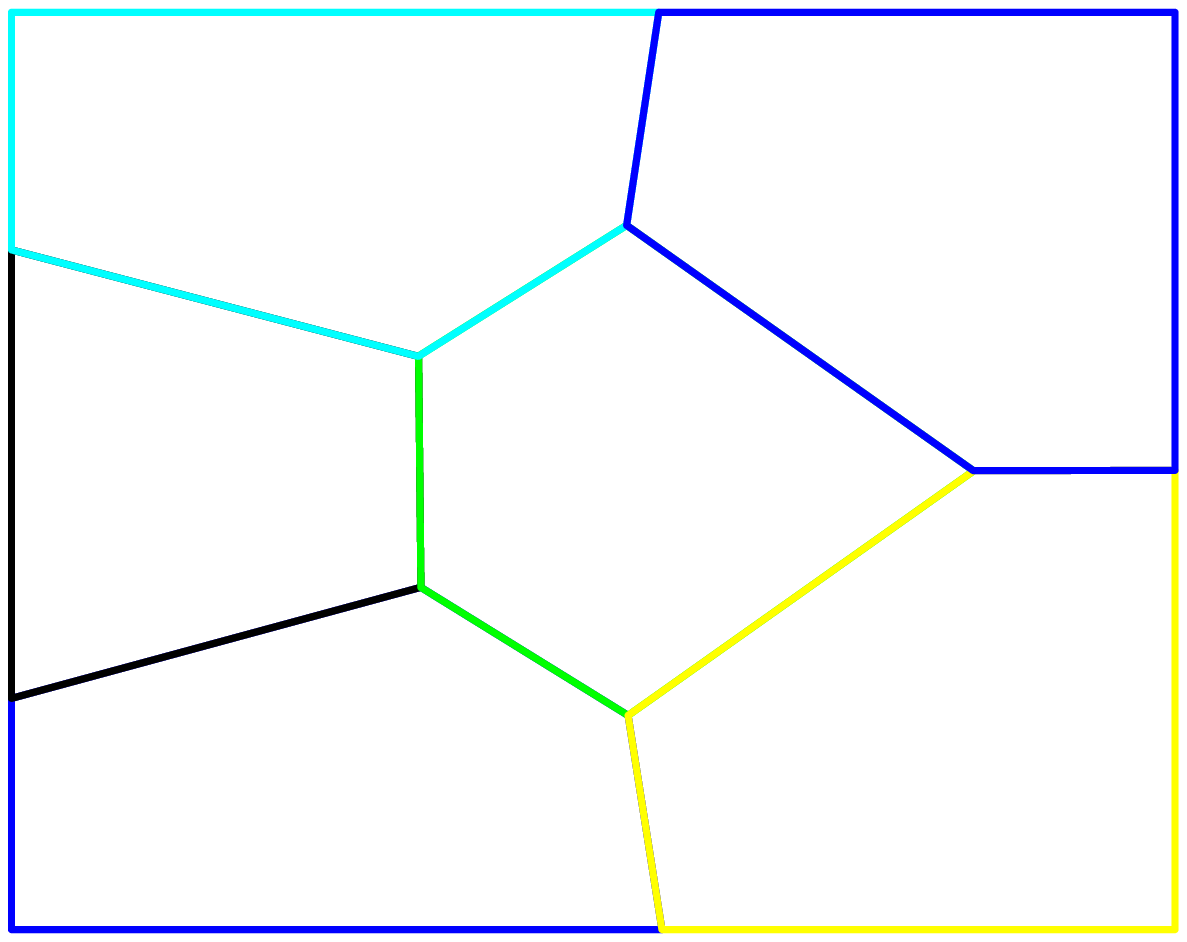}%
\caption{An example simulation of the gossip coverage algorithm with
  uniform random edge selection. The environment $Q$ is a rectangle with
  uniform density, centroids are computed with performance function
  $f(x)=x^2$, and $N=6$ regions compose the partition.  Snapshots of an
  evolution are shown for $t\in\{0,20,50,100,300\}$. One may verify
  numerically that the sequence converges asymptotically to a centroidal
  Voronoi partition.  At $t=20$ one of the regions is disconnected.}
\label{fig:Simulations}
\end{figure*}

We have extensively simulated the partition-based gossip coverage algorithm
$T$ on a $2$-dimensional polygonal environment with uniform density and
performance function $f(x)=x^2$.  Simulations have been implemented as a
{\tt Matlab} program, using the {\tt General Polygon Clipper Library} to
perform operations on polygons.  We adopted the following communication
model: at each iteration, a region pair is chosen, uniformly at random,
among all pairs of adjacent regions.  Fig.~\ref{fig:Simulations} is an
illustration of a typical evolution.

Our first numerical finding is that the gossip coverage algorithm appears
to converge to centroidal Voronoi partitions from all initial
conditions. This is the same property that the Lloyd synchronous coverage
algorithm~\eqref{eq:coverage-control} is known to possess. In other words,
our numerically-computed sequences of partitions always converge to
centroidal Voronoi partitions \textendash{} even though our theoretical
analysis (1) requires a continuous interpolation from $T$ to $\Tdelta$ and
(2) does not exclude convergence to degenerate partitions where some
component regions might have coincident centroids, or might have empty
interiors, or might be composed of ``polygons with an infinite number of
vertices,'' that is, arbitrary sets.


A second numerical finding is that, throughout numerous sample executions,
the resulting polygonal regions rarely have complicated shapes and large
numbers of vertices.  This is good news because of our assumption of finite
convexity and because large numbers of vertices affect both the computation
and the communication burden of the gossip coverage algorithm.



Finally, it is possible, and we have observed it numerically, to have
evolutions of the algorithm that, before converging to centroidal Voronoi
partitions, have components with disconnected regions.  We believe that
there might be applications where it is desirable to maintain connectivity
of the components of the partition and, therefore, we sketch here how to
design a connectivity-preserving algorithm. Note that the update step of
the partition-based coverage algorithm amounts to the exchange, among the
agents, of a region, which consists in general of several connected
components. In the connectivity-preserving algorithm, such components are
considered individually, and each of them is traded only if this can be
done without loosing connectivity; if not, the component is kept by the
previous owner. Numerical simulations indicate that such an algorithm leads
to centroidal Voronoi partitions as well.

\section{On the Krasovskii--LaSalle invariance principle: set-valued maps
  on metric spaces}
\label{sec:Krasovskii-LaSalle}
In this section we consider discrete-time set-valued dynamical system
defined on metric spaces. Our goal is to provide some extensions of the
classical Krasovskii-LaSalle Invariance Principle; we refer the interested
reader to~\cite{AB-LM:05,JH-DL-DA-EDS:05,RGS-RG-ART:08} for recent
invariance principles for switched continuous-time and hybrid systems.

We start by reviewing some preliminary notions including set-valued
dynamical systems, continuity and invariance properties, and Lyapunov
functions.  On a metric space $(X,d)$, where $X$ is a set and $d$ is a
metric on $X$, a set-valued map $\setmap{T}{X}{X}$ is non-empty if
$T(x)\neq \emptyset$ for all $x\in{X}$. An evolution of the dynamical
system determined by a non-empty set-valued map $T$ is a sequence
$\seqdef{x_n}{n\in\integernonnegative}\subset X$ with the property that
\begin{equation*}
  x_{n+1}\in T(x_n), \quad n\in\integernonnegative.
\end{equation*}
In other words, we regard a set-valued map as a nondeterministic
discrete-time dynamical system.  For set-valued maps we introduce notions
of continuity and invariance as follows.  A set-valued map $T$ is
\emph{closed} at $x\in{X}$ if, for all pairs of convergent sequences $x_k
\to x$ and $x'_k \to x'$ such that $x'_k\in T(x_k)$, one has that $x'\in
T(x)$.  Additionally, $T$ is closed on $W\subset X$ if it is closed at all
$w\in{W}$.  A set $W\subset{X}$ is \emph{weakly positively invariant} for
$T$ if $T(w)\intersect{W}$ is non-empty for all $w\in{W}$.  A set $W$ is
\emph{strongly positively invariant} for $T$ if $T(w)\subset{W}$ for all
$w\in{W}$.

We are ready now to state a Krasovskii-LaSalle invariance principle for
set-valued maps defined on metric spaces. This result extends the Global
Convergence Theorem in \cite{DGL:84} to more general Lyapunov functions.
Its proof follows the lines of the proof of Theorem~1.21 in
\cite{FB-JC-SM:09}, and is thus omitted.

\begin{lemma}[Krasovskii-LaSalle invariance principle for set-valued maps]
  \label{lem:Lasalle}
  Let $(X,d)$ be a metric space and $\setmap{T}{X}{X}$ be non-empty. Assume
  that:
  \begin{enumerate}
  \item there exists a compact set $W\subseteq X$ that is strongly positively
    invariant for $T$;
  \item there exists a function $\map{U}{W}{\real}$ such that $U(w') \leq
    U(w)$, for all $w\in W$ and $w'\in T(w)$;
  \item the function $U$ is continuous on $W$ and the map $T$ is closed on $W$.
  \end{enumerate}
  Then there exists $c\in\real$ such that each evolution of $T$ with
  initial condition in $W$ approaches a set of the form $S\intersect
  U^{-1}(c)$, where $S$ is the largest weakly positively invariant set
  contained in
  $$\setdef{w\in{W}} {\exists\, w' \in T(w) \; \text{such that} \;
    U(w')=U(w)}.$$
\end{lemma}

In this paper, given the metric space $(X,d)$, we deal with a set-valued
map $\setmap{T}{X}{X}$ defined by a collection of maps
$\map{T_1,\ldots,T_m}{X}{X}$ via the equality $T(x)=\left\{T_1(x),\ldots,
  T_m(x)\right\}$ for $x\in X$.  For this kind of set-valued maps,
closedness is related to continuity of ordinary maps.

\begin{lemma}\label{lem:Cont->close}
  Let $\map{T_1,\ldots,T_m}{X}{X}$ be continuous on $W\subset{X}$.  The
  set-valued map $\setmap{T}{X}{X}$ defined by $T(x)=\left\{T_1(x),\ldots,
    T_m(x)\right\}$ is closed on $W$.
\end{lemma}
\begin{proof}
  Let $w_n\to w$ and $w'_n\to w'$ be a pair of convergent sequences in $W$,
  such that $w'_n\in T(w_n)$. We claim that the continuity of $T_1,\ldots,
  T_m$ implies $w'\in T(w)$.

  Note that, by hypothesis, for all $n\in \integernonnegative$ there exists
  $i_n\in\left\{1,\ldots,m\right\}$ such that $w'_n=T_{i_n}(w_n)$. Because
  the set $\left\{1,\ldots,m\right\}$ is finite, there exists an index
  $j\in\left\{1,\ldots,m\right\}$ that appears infinitely many times in the
  sequence $\{i_n\}_{n\in\natural}$. Consider the subsequences $\{w_{n_l}\}
  \subseteq \{w_n\}$ and $\{w'_{n_l}\}\subseteq \{w'_n\}$, such that
  $w'_{n_l}=T_{j}(w_{n_l})$. Clearly, we have that $w_{n_l}\to w$ and
  $w'_{n_l}\to w'$, where from the continuity of $T_{j}$ it follows that
  $w'=T_{j}(w)$. Thus, $w'\in T(w)$ and the claim is proved.
\end{proof}

The following result is a stronger version of Lemma~\ref{lem:Lasalle}, for
a particular class of set-valued dynamical systems.

\begin{theorem}[Uniformly persistent switches imply convergence]
  \label{th:StrongerLasalle}
  Let $(X,d)$ be a metric space.  Given a collection of maps
  $\map{T_1,\ldots, T_m}{X}{X}$, define the set-valued map
  $\setmap{T}{X}{X}$ by $T(x)=\left\{T_1(x),\ldots, T_m(x)\right\}$ and let
  $\seqdef{x_n}{n\in\integernonnegative}$ be an evolution of $T$. Assume
  that:
  \begin{enumerate}
  \item\label{item:Compactness} there exists a compact set $W\subseteq X$
    that is strongly positively invariant for $T$;
  \item\label{item:Lyapunov} there exists a function $\map{U}{W}{\real}$ such
    that $U(w')< U(w)$, for all $w\in W$ and $w'\in T(w)\setminus\{w\}$;
  \item\label{item:continuousIterate} the maps $T_i$, for $i\in\until{m}$,
    and $U$ are continuous on $W$; and
  \item\label{item:persistentChoice} for all $i\in \until{m}$, there exists
    an increasing sequence of times $\setdef{n_k}{k\in\integernonnegative}$
    such that $x_{n_k+1}=T_i(x_{n_k})$ and $(n_{k+1}-n_k)$ is bounded.
  \end{enumerate}
  If $x_0\in{W}$, there exists $c\in\real$ such that the evolution
  $\seqdef{x_n}{n\in\integernonnegative}$ approaches the set
  \begin{equation*}
    (F_1\intersection \cdots \intersection F_m) \intersect U^{-1}(c),
  \end{equation*}
  where $F_i=\setdef{w\in W}{T_i(w)=w}$ is the set of fixed points of $T_i$
  in $W$, $i\in\until{m}$.
\end{theorem}

\noindent Loosely speaking, (i) the compactness of a strongly positively
invariant set, (ii) a monotonicity property for a Lyapunov function, (iii)
continuity properties, and (iv) uniformly persistent switches among
finitely many maps, together ensure convergence of each evolution to the
intersection of the fixed points of the maps.

\begin{proof}[Proof of Theorem~\ref{th:StrongerLasalle}]
  Let $S$ be the largest weakly positively invariant set contained in
  \begin{equation*}
    \setdef{w\in W}{\exists\, w'\in T(w) \text{ such that }
      U(w')=U(w)} = F_1\union \cdots \union F_m.
  \end{equation*}
  Since $T$ is closed by Lemma~\ref{lem:Cont->close}, the assumptions of
  Lemma~\ref{lem:Lasalle} are met; hence there exists $c\in\real$ such that
  $U(x_n)\to c$ and $x_n\to S \intersect U^{-1}(c)$.

  Let $\Omega(x_n)$ denote the $\omega$-limit set of the sequence
  $\setdef{x_n}{n\in\integernonnegative}$. If we show that
  $\Omega(x_n)\subseteq (F_1\intersect \cdots \intersect F_m) \intersect
  U^{-1}(c)$, then the statement of the theorem is proved. We proceed by
  contradiction. To this aim, let $\widehat{x}\in S\intersect U^{-1}(c)
  \setminus \big((F_1\intersect \cdots \intersect F_m) \intersect
  U^{-1}(c)\big)$ and let $\seqdef{x_{n_h}}{h\in\integernonnegative}$ be a
  subsequence such that $x_{n_h}\to \widehat{x}$.

  Observe that for each $\widehat{x}\in S \setminus (F_1\intersect \cdots
  \intersect F_m) $, there exists a non-empty set
  $\mathcal{I}_{\widehat{x}} \subset \left\{1,\ldots,m\right\}$ such that,
  $\widehat{x}=T_i\left(\widehat{x}\right)$ if $i\in
  \mathcal{I}_{\widehat{x}}$, and $\widehat{x}\neq
  T_i\left(\widehat{x}\right)$ if $i\notin \mathcal{I}_{\widehat{x}}$. By
  the continuity of the maps $T_i$, there exists $\delta\in\realpositive$
  such that, if $i \notin \mathcal{I}_{\widehat{x}}$, then $T_i(x)\neq x$
  for all $x\in B_\delta(\widehat{x})=\setdef{x\in{W}}{d(x,\widehat{x})\leq
    \delta}$. Let now
 \begin{equation*}
   \gamma_{\delta} = \min_{i\in \mathcal{I}_{\widehat{x}}}
   \left\{\min_{x\in B_{\delta}(\widehat{x})}
     \Big( U\left(x\right)-U\left(T_i(x)\right) \Big)
   \right\} \geq 0.
 \end{equation*}
 By hypothesis, if $i \notin \mathcal{I}_{\widehat{x}}$, then $U(T_i(x)) <
U(x)$ for all $x\in B_\delta(\widehat{x})$.  Hence, since
$B_\delta(\widehat{x})$
 is closed, and $U$ and the maps $T_i$ are continuous, we deduce that
 $\gamma>0$.

 Observe now that hypothesis (iv) implies the existence of a duration
 $D\in \natural$ such that every map $T_i$, $i \in \until{m}$, is applied
 at least once within every interval ${[n,n+D[}$, for all
 $n\in\integernonnegative$.  Consider the set $\left\{T_i\right\}_{i \in
   \mathcal{I}_{\widehat{x}}}$; this is a collection of continuous maps
 having $\widehat{x}$ as fixed point. Then, there exists a suitable $\eps
 \in \realpositive$ such that, given any $r$-uple $(j_1,\ldots,
 j_r)\in\mathcal{I}_{\widehat{x}}^r$, $r\leq D$, we have that
 $T_{j_1}\circ T_{j_2} \circ T_{j_3} \circ \ldots \circ T_{j_r}(x)\in
 B_\delta(\widehat{x})$ for all $w\in B_\eps (\widehat{x})$.

 Select now $k$ such that the element $x_{n_{k}}$ in the subsequence
 $\seqdef{x_{n_h}}{h\in\integernonnegative}$ satisfies $d\left(x_{n_{k}},
   \widehat{x}\right)<\eps$ and ${U(x_{n_k})-c< \gamma_\delta}$.  Let
 \begin{equation*}
   s = \min\setdef{t\in[1,D[}{\exists\, j \notin  \mathcal{I}_{\widehat{x}}
     \text{ such that } x_{n_{k}+t+1}=T_j\left(x_{n_{k}+t}\right)}.
 \end{equation*}
 Observe that $U(x_{n_{k}+s})-c<\gamma_\delta$ and $U(
 x_{n_{k}+s})-U\big(T_j(x_{n_{k}+s})\big)\geq \gamma_\delta$ implying that
 $U\left(T_j\left(x_{n_{k}+s}\right)\right)<c$.  This is a contradiction.
\end{proof}

\begin{remark}
  \label{remark:AndyTeel}
  An alternate proof of this theorem can be given by applying an invariance
  principle obtained in~\cite{RGS-RG-ART:08} on an appropriately-designed
  dynamical extension of the discrete-time set-valued system. \oprocend
\end{remark}

In Appendix~\ref{appA:counterexample} we show how to persistent switching
assumption is necessary.  Next, we provide a probabilistic version of the
previous theorem.

\begin{theorem}[Persistent random switches imply convergence]
  \label{th:StrongerLasalleProb}
  Let $(X,d)$ be a metric space.  Given a collection of maps
  $\map{T_1,\ldots, T_m}{X}{X}$, define the set-valued map
  $\setmap{T}{X}{X}$ by $T(x)=\left\{T_1(x),\ldots, T_m(x)\right\}$.  Given
  a stochastic process $\map{\sigma}{\integernonnegative}{\until{m}}$,
  consider an evolution $\seqdef{x_n}{n\in\integernonnegative}$ of $T$
  satisfying
  \begin{equation*}
    x_{n+1} = T_{\sigma(n)}(x_n).
  \end{equation*}
  Assume that:
  \begin{enumerate}
  \item there exists a compact set $W\subseteq X$ that is strongly positively
    invariant for $T$;
  \item there exists a function $\map{U}{W}{\real}$ such that $U(w')< U(w)$,
    for all $w\in W$ and $w'\in T(w)\setminus\{w\}$;
  \item the maps $T_i$, for $i\in\until{m}$, and $U$ are continuous on $W$;
    and
  \item there exists $p\in{]0,1[}$ and $k\in\natural$ such that, for all
    $i\in\until{m}$ and $n\in\integernonnegative$, there exists
    $h\in\until{k}$ such that
    \begin{equation*}
      \Prob\big[\sigma(n+h)=i \,|\,  \sigma(n),\dots,\sigma(1)\big] \geq p.
    \end{equation*}
  \end{enumerate}
  If $x_0\in{W}$, then there exists $c\in\real$ such that almost surely the
  evolution $\seqdef{x_n}{n\in\integernonnegative}$ approaches the set
  \begin{equation*}
     (F_1\intersection \cdots \intersection F_m) \intersect U^{-1}(c),
  \end{equation*}
  where $F_i=\setdef{w\in W}{T_i(w)=w}$ is the set of fixed points of $T_i$
  in $W$, $i\in\until{m}$.
\end{theorem}

\noindent Loosely speaking, (i) the compactness of a strongly positively
invariant set, (ii) a monotonicity property for a Lyapunov function, (iii)
continuity properties, and (iv) persistent random switches among finitely
many maps, together ensure convergence of each evolution to the
intersection of the fixed points of the maps.

\begin{proof}[Proof of Theorem~\ref{th:StrongerLasalleProb}]
  If $x_0\in W$, then the stochastic process $\sigma$ induces a
  stochastic process taking values in $W$. From now on, we restrict our
  attention to sequences $\seqdef{x_n}{n\in\integernonnegative}$ such that
  $x_0 \in W$. In other words we assume that the sample space containing
  all the evolutions of our interest is given by
 \begin{equation*}
   \mathcal{A}= \bigsetdef{ \seqdef{x_n}{n\in\integernonnegative} }
   {x_n\in W \text{ for all } n \in\integernonnegative }.
 \end{equation*}
 Let $S$ be the largest weakly positively invariant set contained in
 \begin{equation*}
   \setdef{w\in W}{\exists\, w'\in T(w) \text{ such that }
     U(w')=U(w)} = F_1\union \cdots \union F_m.
 \end{equation*}
 From Lemma~\ref{lem:Lasalle}, we know that there exists $c\in\real$ such
 that $x_n\to S\intersect U^{-1}(c)$.  This implies that the
 following event is certain:
 \begin{equation*}
   E = \bigsetdef{\seqdef{x_n}{n\in\integernonnegative}}
   {\exists\, c \in \real \text{ such that } \lim_{n \to \infty} U(x_n)=c }.
 \end{equation*}
 Let $\Omega(x_n)$ denote the $\omega$-limit set of the sequence
$\setdef{x_n}{n\in\integernonnegative}$. If we show that
$\Omega(x_n)\subseteq
 \left((F_1\intersect \cdots \intersect F_m) \intersect U^{-1}(c)\right)$
almost surely, then the statement of the theorem is proved. Next,
consider the event
 \begin{equation*}
   E_1 = \bigsetdef{\seqdef{x_n}{n\in \integernonnegative}}
   {\exists\, \widehat{x}\in S\setminus (F_1\intersect\ldots \intersect F_m)
     \text{ such that } \widehat{x}\in \Omega\left(x_n\right) }.
 \end{equation*}
 Assume by contradiction that $\Prob\left[E_1\right]>0$. Now we compute
 $\Prob\left[E|E_1\right]$.  Note that, for each $\widehat{x}\in S
 \setminus (F_1\intersect \ldots \intersect F_m)$, there exists a
 non-empty set $\mathcal{I}_{\widehat{x}} \subset \until{m}$ such that,
 $\widehat{x}=T_i\left(\widehat{x}\right)$ if $i\in
 \mathcal{I}_{\widehat{x}}$, and $\widehat{x}\neq
 T_i\left(\widehat{x}\right)$ if $i\notin
 \mathcal{I}_{\widehat{x}}$. Similarly to the proof of
 Theorem~\ref{th:StrongerLasalle}, we can associate to each $\widehat{x}$
 a positive real number $\delta$ such that the inequality $x\neq
 T_i\left(x\right)$ holds true for all $x\in
 B_{\delta}(\widehat{x})=\setdef{x\in W} {d(x,\widehat{x})\leq \delta}$
 and for all $i\notin \mathcal{I}_{\widehat{x}}$. Moreover, we can define
 \begin{equation*}
   \gamma_{\delta}=\min_{i\in \mathcal{I}_{\widehat{x}}}
   \left\{\min_{x\in B_{\delta}(\widehat{x})}
     \Big( U(x)-U(T_i(x))\Big) \right\},
 \end{equation*}
 where the continuity of the maps $T_j$, $j\in \until{m}$, and $U$, and
 the closedness of the set $B_{\delta}(\widehat{x})$ ensure that
 $\gamma_{\delta}>0$.

 Consider the set $\left\{T_i\right\}_{i \in \mathcal{I}_{\widehat{x}}}$;
 this is a collection of continuous maps having $\widehat{x}$ as fixed
 point. Therefore, there exists a suitable $\eps \in \realpositive$ such
 that, given any $r$-uple $(j_1,\ldots, j_r)\in\mathcal{I}_{\widehat{x}}
 ^r$, $r\leq k$, we have $T_{j_1}\circ T_{j_2} \circ T_{j_3} \circ \ldots
 \circ T_{j_r}(x)\in B_\delta(\widehat{x})$ for all $x\in B_\eps
 (\widehat{x})$.  Given $\seqdef{x_n}{n\in\integernonnegative}$, if there
 exists $\widehat{x}\in S\setminus (F_1\intersect\ldots \intersect F_m)$
 such that $\widehat{x}\in \Omega\left(x_n\right)$, then there must exist
 $\{n_h|\; h \in \integernonnegative\}$ such that $x_{n_h}\in
 B_{\eps}(\hat{x})$ for all $h\in \integernonnegative$. Moreover, without
 loss of generality we can assume that $n_{h+1}-n_h>k$ for all $h \in
 \integernonnegative$.  Consider now the event
 \begin{multline*}
   E_3 = \big\{ \setdef{i_n\in\until{m}}{n\in\integernonnegative} \;|
   \\
   \exists\, \bar{h} \text{ such that } i_{n_h+s}\in I_{\widehat{x}}
     \text{ for all } s\in\until{k} \text{ and } h\geq
     \bar{h} \big\}.
 \end{multline*}
 To compute $\Prob[E_3]$, we define, for $j\in \integernonnegative$,
 \begin{equation*}
   E_{3,j} = \bigsetdef{ \setdef{i_n\in\until{m}}{n\in\integernonnegative}}
   {i_{n_h+s}\in I_{\widehat{x}} \text{ for all }
     s\in\until{k} \text{ and } h\geq j }.
 \end{equation*}
 Observe that $\setdef{E_{3,j}}{j \in \integernonnegative}$ is a countable
 collection of disjoint sets such that $E_3=\bigcup_{j=0}^\infty
 E_{3,j}$. By hypothesis we have that
 \begin{equation*}
   \Prob[E_{3,j}]\leq\lim_{l \to \infty}\prod_{s=j}^{l}(1-p)=0,
 \end{equation*}
 and therefore $\Prob[E_3]=0$. This implies that, almost surely, there
 exists a subsequence $\seqdef{n_{h_s}}{s\in\integernonnegative} \subseteq
 \seqdef{n_h}{h\in\integernonnegative}$ with the property that, for all
 $s\in\integernonnegative$, $x_{n_{h_s}+1}=T_i(x_{n_{h_s}})$ for some $i
 \notin \mathcal{I}_{\widehat{x}}$ and, therefore, also with the property
 that $U(x_{n_{h_s}})- U(x_{n_{h_k}+1})> \gamma_{\delta}$. Consequently,
 almost surely, we have that $\lim_{s\to \infty} U(x_{n_{h_s}})=-\infty$
 thus violating the fact that $E$ is a certain event. This implies that
 $\Prob\left[E_1\right]=0$ and that, almost surely, $x_n \to
 (F_1\intersect \ldots \intersect F_m) \intersect U^{-1}(c)$.
\end{proof}

\begin{remark}\label{rem:NoStrongly}
  The assumption, in Lemma~\ref{lem:Lasalle} and
  Theorems~\ref{th:StrongerLasalle} and~\ref{th:StrongerLasalleProb}, that
  $W$ is strongly positively invariant ensures that any evolution of $T$
  with initial condition in $W$ remains in $W$. By relaxing this
  assumption, it is possible to provide weaker versions of these
  results. Specifically, requiring $W$ to be only compact and not
  necessarily strongly positively invariant, the thesis of
  Lemma~\ref{lem:Lasalle} and Theorems~\ref{th:StrongerLasalle}
  and~\ref{th:StrongerLasalleProb} do not hold in general for any evolution
  of $T$ with initial condition on $W$, but are still valid for those
  evolutions $\seqdef{x_n}{n\in\integernonnegative}$ that take values in $W$
  for all $n\in\integernonnegative$. \oprocend
\end{remark}

\section{On the topology of the space of partitions: compactness properties
  in the symmetric difference metric}\label{sec:PartitionsSpace}

Motivated by the invariance principles presented in
Section~\ref{sec:Krasovskii-LaSalle}, we study metric structures on the set
of partitions, with a focus on compactness and continuity
properties. Specifically, we show how a particular subset of the set of
partitions can be regarded as a compact metric space and how certain
relevant maps are continuous over that subspace.  In this section, and only
in this section, the assumptions on $Q$ are relaxed to give more general
results: we assume that $Q\subset\real^d$ is compact and connected and has
non-empty interior.

Let $\C$ denote the {\em set of closed subsets} of $Q$.  We would like to
introduce a topology on $\C$ with two properties: $\C$ is compact and the
Voronoi map, the centroid map, and the multicenter function, defined in
equations~\eqref{eq:DefVor} \eqref{eq:def-Cd}, and \eqref{eq:H-centroid}
respectively, are continuous over $\C$. A natural candidate is the topology
induced by the well-known~\cite{GS-RJBW:79} Hausdorff metric on $\C$: given
two sets $A, B \in\C$, their Hausdorff distance is
$d_H(A,B)=\max\left\{\max_{a\in A}\min_{b\in B} d(a,b), \max_{b\in B}
  \min_{a\in A} d(a,b)\right\}$.
This metric induces a topology on $\C$ which makes it a compact space, but
is not suitable for our purpose because, with respect to this topology, the
Voronoi map, the centroid map, and the multicenter function, are not
continuous; see Appendix~\ref{app:discontinuous}. Additionally, note that,
unlike the Hausdorff metric, the centroid map and the multicenter function
are insensitive to sets of measure zero.

In what follows, we introduce the symmetric difference metric, as a metric
insensitive to sets of measure zero.  Given two subsets $A,B\in\C$, we
define their \emph{symmetric difference} by $A\ominus B=(A\cup B)\setminus
(A\intersect B)$. Moreover, letting $\measure$ denote the Lebesgue measure
on $\real^d$, we define the \emph{symmetric difference distance}, also
called the \emph{symmetric distance} for simplicity,
$\map{d_\ominus}{\C\times\C}{\realnonnegative}$ by
\begin{equation*}
  d_\ominus(A,B)=\measure(A\ominus B),
\end{equation*}
that is, the symmetric distance between two sets is the measure of their
symmetric difference. Given these notions, it is useful to identify sets
that differ by a set of measure zero: we write $A\sim B$ whenever
$\measure(A\ominus{B})=0$. Clearly, $\sim$ is an equivalence relationship
on $\C$ and, accordingly, we let $\C^*=\C/\!\sim$ denote the \emph{quotient
  set of closed subsets of $Q$}.  Now, for any two elements $A^*$ and $B^*$
of $\C^*$, we define $d_\ominus(A^*, B^*)=d_\ominus(A,B)$ where $A$ and $B$
are any representatives of $A^*$ and $B^*$, respectively. With this notion
of $d_\ominus$ on $\C^*\times \C^*$, it is easy to verify that $(\C^*,
d_\ominus)$ is a metric space.  However, to the best of our knowledge no
compactness result is available for $(\C^*, d_\ominus)$.

Next, we introduce a particular family of subsets of $\C$ whose structure
is sufficiently rich for our algorithm.
For $\ell\in\natural$, let $\CL\subset\C$ denote the set of
\emph{$\ell$-convex and closed subsets} of $Q$, that is, the set of subsets
of $Q$ equal to the union of $\ell$ convex and closed subsets of
$Q$. Formally, we set
\begin{equation*}
  \CL = \bigsetdef{v\subseteq Q }{ v=\union_{i=1}^\ell S_i \text{ where
    } S_1,\ldots,S_{\ell} \subseteq Q \text{ are convex and closed} }.
\end{equation*}
Note that we do not require the sets $S_1,\dots,S_\ell$ to be distinct so
that $\C_{(k)} \subset \CL$, for any $k<\ell$.  In what follows we study
the \emph{quotient set of $\ell$-convex and closed subsets}
$\CL^*=\CL/\!\sim$. The next result is the main result of this section.

\begin{theorem}[Compactness of $\CL^*$]\label{th:CL}
  The pair $(\CL^*,d_\ominus)$ is a metric space and, with the topology
  induced by $d_\ominus$, the set $\CL^*$ is compact.
\end{theorem}
\begin{proof}
  It is easy to verify that $d_\ominus$ is a metric on $\C^*$.  Instead,
  proving the compactness of $\CL^*$ requires some attention.  We aim to
  show that any sequence in $\CL^*$ has a subsequence that converges to a
  point in $\CL^*$. This fact's proof is articulated in several steps and
  relies upon several known results:
  \begin{enumerate}

  \item \label{claim:C-is-H-compact} the space $\C$, endowed with the
    Hausdorff distance $\map{d_H}{\C\times\C}{\realnonnegative}$,
    is~\cite[Theorem~0.8]{SBN:78} a compact metric space;
%
%

  \item \label{claim:Hconv->convex} if a sequence of closed convex subsets
    of $Q$ converges in the Hausdorff sense to a set $K$,
    then~\cite[Proposition~1.6.8]{SGK-HRP:08} $K$ is closed and convex; and

  \item \label{claim:Delta-H} for any two convex subsets $A, B$ of
    $Q\subset\real^d$, it is known~\cite[Eq.~(1)]{HG:00} that
    $$ d_\ominus(A,B)\leq
    \Big(\frac{2\kappa_d}{2^{1/d}-1}\left(\frac{D}{2}\right)^{d-1}
    \Big)d_H(A,B),$$
    where $D=\max\left\{\diam(A),\diam(B)\right\}$ and where $\kappa_d$ is
    the volume of the unit sphere in $\real^d$.
\end{enumerate}
Now let $\seqdef{v^*(n)}{n\in\integernonnegative}$ be any sequence in
$\CL^*$. For all $n\in\integernonnegative$, pick any representative of the
equivalence class $v^*(n)$, denote it by $v(n)$ and consider the sequence
$\seqdef{v(n)}{n\in\integernonnegative}$. Because
$\seqdef{v(n)}{n\in\integernonnegative}$ is a sequence in $\CL$ and because
$\CL\subseteq \C$, it follows from fact~(i) above that
$\seqdef{v(n)}{n\in\integernonnegative}$ contains a subsequence that
converges in the Hausdorff sense to an element of $\C$. In other words,
there exist $\seqdef{v(n_k)}{k\in\integernonnegative} \subseteq
\seqdef{v(n)}{n\in\integernonnegative}$ and $\bar{v}\in\C$ such that
$$
\lim_{k\to \infty} v(n_k)\stackrel{\text{{\tiny (H)}}}{=}\bar{v},
$$
where $\stackrel{\text{{\tiny (H)}}}{=}$ denotes convergence in the
Hausdorff sense.  We claim that $\bar{v}\in\CL$ so that the set $\CL^*$ is
compact. To show this claim, we plan to exhibit a collection of convex and
closed subsets of $Q$, say $\left\{S_1,\ldots,S_{\ell}\right\}$, such that
$\bar{v}= \cup_{i=1}^{\ell} S_i$.

We begin to prove this claim as follows.  By definition of $\CL$, for all
$k\in \integernonnegative$ there exists a collection $\{S_{1}(n_k), \ldots,
S_{\ell}(n_k)\}$, of convex and closed subsets of $Q$ whose union equals
$v(n_k)$.  Now, we consider the sequence
$\seqdef{S_{1}(n_k)}{k\in\integernonnegative}$. Again, since $(\C, d_H)$ is
a compact metric space we have that there exists a subsequence
$\seqdef{n_{k_1}}{k_1\in\integernonnegative} \subseteq
\seqdef{n_{k}}{k\in\integernonnegative}$ such that $\lim_{k_1\to\infty}
S_1(n_{k_1}) \Hequal \bar{S}_1$ for some
$\bar{S}_1$. Fact~\ref{claim:Hconv->convex} above ensure that $\bar{S}_1$
is a convex closed subset of $Q$. Consider now the sequence
$\seqdef{S_2(n_{k_1})}{k_1\in\integernonnegative}$. By reasoning as
previously we know that there exists a subsequence
$\seqdef{n_{k_2}}{k_2\in\integernonnegative} \subseteq
\seqdef{n_{k_1}}{k_1\in\integernonnegative}$ such that $\lim_{k_2\to\infty}
S_2(n_{k_2}) \Hequal \bar{S}_2$ where $\bar{S}_2$ is some convex closed set
of $Q$. Moreover, it is clear that also $\lim_{k_2\to\infty} S_1(n_{k_2})
\Hequal \bar{S}_1$.  By iterating this procedure we conclude that there
exist a sequence $\seqdef{n_s}{s\in\integernonnegative}$ and a collection
of convex and closed sets $\{\bar{S}_1,\ldots,\bar{S}_\ell\}$ such that
$\lim_{s\to\infty} S_i(n_s) \Hequal \bar{S}_{i}$ for all
$i\in\until{\ell}$.

Next, we aim to show that
\begin{equation}\label{eq:InvLimUnion}
  \lim_{s\to\infty} \bigcup_{i=1}^\ell S_i(n_s) %
  = \bigcup_{i=1}^{\ell}  \lim_{s\to\infty} S_i(n_s) %
  =\bigcup_{i=1}^{\ell}\bar{S}_i.
\end{equation}
For simplicity, let us denote
$\lim_{s\to\infty}\bigcup_{i=1}^\ell S_i(n_s)$ by $S_{\infty}$. For $p\in
Q$, note
\begin{equation*}\label{eq:Dist_p}
  p\in S_\infty
  \iff
 \dist\left(p, \lim_{s\to\infty} \bigcup_{i=1}^\ell S_i(n_s)\right)=0
  \iff
  \lim_{s\to\infty} \dist\left(p, \bigcup_{i=1}^\ell S_i(n_s)\right)=0,
\end{equation*}
where, for a given closed set $X$, $\dist(p,X)$ denotes the Euclidean
distance between $p$ and $X$, namely, $\min_{x\in X}\| p-x\|$. Using the
fact that, for given closed sets $X$ and $Y$, $\dist \left(p, X \cup
  Y\right)=\min \left\{\dist(p, X), \dist(p, Y)\right\}$, we can write
\begin{align*}
  \lim_{s\to \infty} \dist \left(p, \bigcup_{i=1}^\ell S_i(n_s)\right)=0
  &\enspace\iff\enspace \lim_{s\to \infty} \min\left\{\dist(p, S_1(n_s)),
    \ldots, \dist(p, S_\ell(n_s)) \right\}=0\\
  &\enspace\iff\enspace \exists j\in \until{\ell} \text{ s.t. } \lim_{s\to \infty} \dist(p, S_j(n_s))=0\\
  &\enspace\iff\enspace  \exists j\in \until{\ell} \text{ s.t. }  p\in \bar{S}_j\\
  &\enspace\iff\enspace p \in \bigcup_{i=1}^\ell \bar{S}_i
\end{align*}
The above chain of implications proves~\eqref{eq:InvLimUnion}. Now observe
that, from the uniqueness of the limit it follows that $\lim_{k \to \infty}
v(n_k)=\lim_{s\to \infty} \bigcup_{i=1}^\ell S_i(n_s)=\bigcup_{i=1}^\ell
\bar{S}_i=\bar{v}$.  Since $\bar{S}_i$ is closed and convex for all $i\in
\until{\ell}$, it follows that $\bar{v}\in \CL$ and, in turn, that $\CL$
endowed with the Hausdorff metric is a metric compact space.

To establish the statement of the Theorem it finally remains to prove that
either $\lim_{k\to \infty} d_\ominus(v(n_k), \bar{v})=0$ or $\lim_{s\to
  \infty} d_\ominus(v(n_s), \bar{v})=0$. To this end, observe that, given
$X_1, X_2, Y_1, Y_2 \subseteq Q$, the following inclusion holds $(X_1\cup
X_2) \ominus (Y_1\cup Y_2) \subseteq (X_1 \ominus Y_1) \cup (X_2 \ominus
Y_2)$. Hence, we compute
\begin{align*}
  v\left(n_s\right)\ominus \bar{v}=\left(\bigcup_{i=1}^\ell S_i(n_s)\right)
  \ominus \bar{v} = \left(\bigcup_{i=1}^\ell S_i(n_s)\right) \ominus
  \left(\bigcup_{i=1}^\ell \bar{S}_i\right) \subseteq \bigcup_{i=1}^{\ell}
  \left(S_i(n_s) \ominus \bar{S}_i\right),
\end{align*}
which implies
\begin{align*}
  d_\ominus\left(v\left(n_s\right),\bar{v}\right) %
  &=d_\ominus \left(\bigcup_{i=1}^\ell S_i(n_s), \bar{v}\right)
  \leq \sum _{i=1}^\ell d_\ominus \left(S_i(n_s), \bar{S}_i\right)\\
  & \leq \sum _{i=1}^\ell
  \frac{2\kappa_d}{2^{1/d}-1}\left(\frac{\max\left\{\diam\left(S_i(n_s)\right),
        \diam\left(\bar{S}_i\right)\right\}}{2}\right)^{d-1} d_H
  \left(S_i(n_s), \bar{S}_i\right)\\
  &\leq \frac{2\kappa_d}{2^{1/d}-1}\left(\frac{\diam(Q)}{2}\right)^{d-1}
  \sum _{i=1}^\ell    d_H \left(S_i(n_s), \bar{S}_i\right),
\end{align*}
where the second and the third inequalities follow, respectively, from
fact~\ref{claim:Delta-H} above and from the upper bounds
$\max\{\diam\left(S_i(n_s)\right), \diam\left(\bar{S}_i\right)\}\leq
\diam(Q)$, for all $i\in \until{\ell}$. Since $\lim_{s\to \infty} d_H
\left(S_i(n_s), \bar{S}_i\right)=0$ for all $i\in \until{\ell}$, we
conclude that
\begin{equation*}
  \lim_{s\to \infty} d_\ominus(v(n_s), \bar{v})=0.
\end{equation*}
Now, let $\bar{v}^*$ denote the equivalence class for which $\bar{v}$ is a
representative. Since the metric $d_\ominus$ is insensitive to sets of
measure zero, it follows that $\lim_{s\to \infty} d_\ominus(v^*(n_s),
\bar{v}^*)=0$ and, in turn, that $\seqdef{v^*(n)}{n\in\integernonnegative}$
has a subsequence convergent to point of $\CL^*$. This concludes the proof
that $\CL^*$ is a compact space.
\end{proof}

The metric $d_\ominus$ naturally extends to a metric over the space
$(\C^*)^N$, and hence over $(\CL^*)^N$, by defining
\begin{equation}\label{eq:d}
  d_\ominus(u,v)=\sum_{i=1}^{N}d_\ominus(u_i,v_i),
\end{equation}
for any $u=(u_i)_{i=1}^N$ and $v=(v_i)_{i=1}^N$ in $(\C^*)^N$.  The
compactness of $(\CL^*)^N$ is then a simple consequence of
Theorem~\ref{th:CL}.

\begin{corollary}[Compactness of $(\CL^*)^N$]
  \label{corol:PartSpaceCompact}
  The pair $\big((\CL^*)^N,d_\ominus\big)$ is a metric space and, with the
  topology induced by $d_\ominus$, $(\CL^*)^N$ is compact.
\end{corollary}

\section{On the continuity of some geometric maps and the resulting
  convergence proofs}
\label{sec:ContinuityMaps}
In this section we prove the main convergence theorems for our gossip
coverage algorithms.  First, however, we need to establish the continuity
properties of certain geometric maps and of the proposed algorithms $T$ and
$\Tdelta$.  Before proceeding, we discuss two significant modeling
aspects. First, recall that Section~\ref{sec:PartitionsSpace} introduces
the spaces $\C$, $\CL$, $\C^N$, $\CL^N$ and the corresponding quotient sets
$\C^*$, $\CL^*$, $(\C^*)^N$, $(\CL^*)^N$. We can do the same with the
partition space $\part$. Indeed from
Definition~\ref{def:partition}\ref{def:Non-EmptyPart} we have that each
component $v_i$ of $v\in \part$ can be mapped by the canonical projection
into an equivalence class $v_i^*$ in $\C^*\setminus\{\emptyset\}$; in turn,
any $v \in \part$ can be naturally associated to the $N$-collection of
equivalence classes $v^*=(v_i^*)_{i=1}^N$ . Accordingly, we denote the
\emph{space of equivalence classes of $N$-partitions} by $\part^*$. Recall
that all these quotient spaces are metric spaces with respect to the
symmetric distance $d_\ominus$.

Second, recall the following maps: the centroid map $\map{\Cd}{
  \setdef{A\in\C}{\mu(A)>0} }{Q}$ defined in equation~\eqref{eq:def-Cd},
the $1$-center function $\map{\Hone}{Q\times\C}{\realnonnegative}$ defined
in equation~\eqref{eq:H1def}, the multicenter function
$\map{\Hcentroid}{\part}{\realnonnegative}$ defined in
equation~\eqref{eq:H-centroid}, the gossip coverage map
$\map{T_{ij}}{\part}{\part}$ with $i\not=j$ defined in
equation~\eqref{eq:OurAlgo-ij}, and, for $\delta>0$ and $i\not=j$, the
modified gossip coverage map $\map{\Tdelta_{ij}}{\part}{\part}$ defined in
Section~\ref{sec:SmoothedMap}.  We claim that all these maps are
insensitive to sets of measure zero.  To substantiate this claim, observe
that the integrals of a bounded measurable function over a set $A$ and over
a set $B$ are equal if $d_\ominus(A,B)=0$.  This observation allows us to
redefine the centroid map, the $1$-center function and the multicenter
function as $\map{\Cd}{\C^*\setminus\{\emptyset\}}{Q}$,
$\map{\Hone}{Q\times\C^*}{\realnonnegative}$ and
$\map{\Hcentroid}{\part^*}{\realnonnegative}$, respectively.  Regarding the
modified gossip coverage map $\Tdelta_{ij}$, we reason as follows.  For
$v^*\in\part^*$, let $v\in\part$ and $v'\in\part$ be two representatives of
$v^*$ and let $\widehat{v}$ and $\widehat{v}'$ denote, respectively, the
images of $v$ and $v'$ under the map $\Tdelta_{ij}$, that is,
$\widehat{v}=\Tdelta_{ij}(v)$ and $\widehat{v}'=\Tdelta_{ij}(v')$. Since
the centroid map and the definitions of the points $\hat{p}_i$ and
$\hat{p}_j$ in equation~\eqref{eq:pi-pj} are insensitive to sets of measure
zero, it follows that $d_\ominus(\widehat{v},\widehat{v}')=0$; in other
words $\widehat{v}$ and $\widehat{v}'$ belong to the same equivalence
class, say $\widehat{v}^*\in \part^*$.  From these facts we can redefine
the modified gossip coverage map as $\map{\Tdelta_{ij}}{\part^*}{\part^*}$.
An analogous argument applies to the map~$T$.  This concludes the
justification of our claim.  Finally, note that the Voronoi map
$\map{V}{Q^N\setminus{S_N}}{\part}$ defined in equation~\eqref{eq:DefVor}
can be composed with the standard quotient projection map and therefore
denoted by $\map{V}{Q^N\setminus{S_N}}{\part^*}$.  In summary, the centroid
map, the $1$-center function, the multicenter function, the modified gossip
coverage map, and the Voronoi map are indeed insensitive to sets of measure
zero.

\subsection{Continuity of various geometric maps}

We start by recalling that the compact connected set $Q$ is equipped with a
bounded measurable positive function $\map{\phi}{Q}{\realpositive}$. We
define the diameter of $Q$ and the infinity norm of $\phi$ by
$\diam(Q)=\max\setdef{\|x-y\|}{x,y\in{Q}}$ and $\norminf{\phi} = \max_{x\in
  Q}{\phi(x)}$, respectively.  The following lemma states some important
properties of the $1$-center cost function.

\begin{lemma}[Continuity properties of the $1$-center
  function]\label{lem:H1cont}
  Given a compact convex set $Q\subset \real^d$, let
  $\map{\phi}{Q}{\realpositive}$ be bounded and measurable and let
  $\map{f}{\realnonnegative}{\realnonnegative}$ be locally Lipschitz,
  increasing, and convex.  Define the function
  $\map{\Hone}{Q\times\C^*}{\realnonnegative}$ as in
  equation~\eqref{eq:H1def}.  The following statements hold:
  \begin{enumerate}
  \item\label{item:p-Hconv} the function $p\mapsto \Hone(p;A)$ is strictly
    convex in $p$, for any $A\in\C^*\setminus\{\emptyset\}$,
  \item\label{item:p-Hcont} the function $p\mapsto \Hone(p;A)$ is globally
    Lipschitz in $p$, for any $A\in\C^*$, and
  \item\label{item:A-Hcont} the function $A\mapsto \Hone(p;A)$ is globally
    Lipschitz in $A$ with respect to $d_\ominus$, for any $p\in{Q}$.
  \end{enumerate}
\end{lemma}

\noindent We now state the main result of this subsection.

\begin{theorem}[Continuity properties of centroid, multicenter and Voronoi
  maps]\label{th:EveryThingCont}
  Given a compact convex set $Q\subset \real^d$, let
  $\map{\phi}{Q}{\realpositive}$ be bounded and measurable and let
  $\map{f}{\realnonnegative}{\realnonnegative}$ be locally Lipschitz,
  increasing, and convex. With respect to the topology induced by
  $d_\ominus$, the following maps are continuous:
  \begin{enumerate}
  \item\label{item:GenCentroidCont} the centroid map
    $\map{\Cd}{\C^*\setminus\{\emptyset\}}{Q}$,
  \item\label{item:Hcentroid-is-cont} the multicenter function
    $\map{\Hcentroid}{\part^*}{\realnonnegative}$,
  \item\label{item:V-cont} the Voronoi map $\map{V}{Q^N\setminus
      S_N}{\part^*}$,
  \item\label{item:TCont} for all $i,j\in \until{N}$, $i\not=j$, the gossip
    coverage map $\map{T_{ij}}{ \setdef{v\in\part^*}{\Cd(v_i)\neq\Cd(v_j)}
    }{\part^*}$, and
  \item\label{item:TepsCont} for all $\delta>0$, $i,j\in \until{N}$,
    $i\not=j$, the modified gossip coverage map
    $\map{\Tdelta_{ij}}{\part^*}{\part^*}$.
  \end{enumerate}
\end{theorem}
The continuity properties~\ref{item:Hcentroid-is-cont} and~\ref{item:TCont}
(respectively,~\ref{item:TepsCont}) are exactly what is needed to apply the
Krasovskii-LaSalle invariance principles stated in
Section~\ref{sec:Krasovskii-LaSalle} to the gossip coverage algorithm
(respectively, to the modified gossip coverage algorithm). The continuity
properties~\ref{item:GenCentroidCont} and~\ref{item:V-cont} are
intermediate results of independent interest.

\begin{proof}[Proof of Lemma~\ref{lem:H1cont}]
  Let $L_f$ be the Lipschitz constant of
  $\map{f}{[0,\diam(Q)]}{\realnonnegative}$.  We check the claims in order,
  beginning with statement~\ref{item:p-Hconv}. For $\lambda\in (0,1)$,
  $A\in\C^*\setminus\{\emptyset\}$, and $p_1,p_2\in Q$, we compute
  \begin{align}
    \Hone(\lambda p_1+ (1-\lambda) p_2;A)
    &=\int_{A}f(\|\lambda p_1+ (1-\lambda) p_2-q\|)\phi(q)dq \nonumber \\
    &\le\int_{A}f(\lambda\|p_1-q\|+ (1-\lambda) \|p_2-q\|)\phi(q)dq
    \label{ineq:triangle+increasing} \\
    &\le \int_{A}\Big(\lambda f(\|p_1-q\|)+ (1-\lambda)
    f(\|p_2-q\|)\Big)\phi(q)dq \label{ineq:convex}\\
    &=\lambda \Hone(p_1;A)+ (1-\lambda) \Hone(p_2;A),\nonumber
  \end{align}
  where inequality~\eqref{ineq:triangle+increasing} follows from the
  triangle inequality and from $f$ being increasing, and
  inequality~\eqref{ineq:convex} follows from the convexity of $f$.
  This inequality
  proves convexity. Moreover, since the first inequality is strict outside
  the line passing through $p_1$ and $p_2$ and since $A$ has non-empty
  interior, the function is in fact strictly convex. Note that
  statement~\ref{item:p-Hconv} implies that $p\mapsto \Hone(p;A)$ is
  locally Lipschitz, using \cite[Theorem~10.4]{RTR:70}. The stronger
  statement~\ref{item:p-Hcont} can be derived as follows. For $p_1,p_2\in
  Q$, we compute
  \begin{align*}
    |\Hone(p_1;A)-\Hone(p_2;A)|
    &=\left|\int_{A}f(\|p_1-q\|)\phi(q)dq - \int_{A}f(\|p_2-q\|)\phi(q)dq\right|\\
    &=\left|\int_{A}[f(\|p_1-q\|) - f(\|p_2-q\|)]\phi(q)dq\right|\\
    &\le\int_{A}|f(\|p_1-q\|) - f(\|p_2-q\|)|\phi(q)dq\\
    &\le\int_{A} L_f \|p_1-p_2\|\phi(q)dq \;\le\; L_f \norminf{\phi} \measure(A) \|p_1-p_2\|.
  \end{align*}
  This implies the Lipschitz condition in statement~\ref{item:p-Hcont}.
  Statement~\ref{item:A-Hcont} can be proved as follows. Let $A$, $A'$ be
  two elements of $\C^*$, note $A=(A\setminus A')\union(A\intersect A')$
  and compute
  \begin{align*}
    \left|\Hone(p;A)-\Hone(p;A')\right|
    &=\left|\int_{A\setminus A'}f(\|p-q\|)\phi(q)dq - \int_{A'\setminus
        A}f(\|p-q\|)\phi(q)dq\right|\\
    &\le\left|\int_{A\setminus A'}f(\|p-q\|)\phi(q)dq\right|+\left|
      \int_{A'\setminus A}f(\|p-q\|)\phi(q)dq\right|\\
    &=\int_{A\setminus A'}f(\|p-q\|)\phi(q)dq+ \int_{A'\setminus
      A}f(\|p-q\|)\phi(q)dq\\
    &=\int_{A\ominus A'}f(\|p-q\|)\phi(q)dq\\
    &\le  \max\setdef{f(\|p-q\|)}{p,q\in A\ominus A'} \, \norminf{\phi} \,
    \measure{(A\ominus A')}\\
    & \le f(\diam(Q)) \, \norminf{\phi} \, d_\ominus(A,A'),
  \end{align*}
  where last inequality follows from $f$ being increasing.  The bound
  implies the Lipschitz condition.
\end{proof}

Before proving Theorem~\ref{th:EveryThingCont} we need the following lemma
about perturbations of convex optimization problems.
\begin{lemma}\label{lem:argminCont}
  Given a compact convex set $Q\subset\real^d$ and a metric space $(X,d)$,
  let $\map{H}{Q\times{X}}{\real}$ have the properties that
  \begin{enumerate}
  \item the map $x\mapsto H(q,x)$ is globally Lipschitz for all $q\in Q$,
    and
  \item the map $q\mapsto H(q,x)$ is continuous and strictly convex.
  \end{enumerate}
  Then the map $\map{q^*}{X}{Q}$, defined by $q^*(x)=\argmin_{q\in
    Q}{H(q,x)}$, is continuous.
\end{lemma}
\begin{proof}
  Let $L_H$ be the Lipschitz constant of $x\mapsto H(q,x)$.  Thanks to the
  Lipschitz condition of the function $x\mapsto H(q,x)$, for all $x,y\in
  X$, the point of minimum $q^*(x)$ takes value in $S=\setdef{q\in
    Q}{H(q,y)\le H(q^*(y),y)+2 L_H d(y,x)}$. Since $S$ is a sub-level set
  of the strictly convex function $q\mapsto H(q,y)$, and since the diameter
  of a sub-level set depends continuously on the level, the distance
  $\|q^*(y)-q^*(x)\|$ can be made arbitrary small by reducing
  $d(y,x)$. This implies the claimed continuity.
\end{proof}

\noindent We are now ready to prove the main result.
\begin{proof}[Proof of Theorem~\ref{th:EveryThingCont}]
  We prove the theorem claims in the order in which they are presented.
  Claim~\ref{item:GenCentroidCont} follows combining
  Lemma~\ref{lem:argminCont} and Lemma~\ref{lem:H1cont}.  Since the
  multicenter function $\Hcentroid$ is a sum of suitable $1$-center
  functions $\Hone$, the claim~\ref{item:Hcentroid-is-cont} is also
  immediate.

  Regarding claim~\ref{item:V-cont}, we discuss in detail the two
  dimensional case. We first prove the continuity when $N=2$: let $p_1$ and
  $p_2$ be points in $Q$. Let $2 l=\|p_1-p_2\|$. Since $p_1\neq p_2$,
  $l>0$. Up to isometries, we can assume that, in the Euclidean plane
  $(x,y)$, $p_1=(-l,0)$ and $p_2=(l,0)$. Let $d_1$ and $d_2$ be the
  distances from the origin of points $p_1$ and $p_2$, respectively. It is
  clear that the two Voronoi regions of $p_1$ and $p_2$ are separated by
  the locus of points $\setdef{x\in Q}{\|x-p_1\|=\|x-p_2\|}$, that is the
  vertical axis.  Now, we assume that the positions of $p_1$ and $p_2$ are
  perturbed by a quantity less than or equal to $\eps$, with
  $0<\eps<l$. By effect of the perturbation, the axis separating the two
  Voronoi regions is perturbed, but it is contained in the locus of points
  $Y_{1 2}(\eps)=\setdef{x\in Q}{ \|x-p_1\| - \|x-p_2\| \le 2\eps}$. By
  definition, this is the set comprised between the two branches of the
  hyperbola whose equation is
  $\frac{x^2}{\eps^2}-\frac{y^2}{l^2-\eps^2}=1$. By elementary
  geometric considerations, the area of this region can be upper bounded by
  \begin{align*}
    \measure(Y_{1 2}(\eps))&\leq 2\, \eps\, 2\, \diam(Q)+ 4 \diam(Q)^2
    \frac{{\eps}/{l}}{\sqrt{1-\frac{\eps^2}{l^2}}} \le 4\,
    \diam(Q)\Big(1+\frac{\diam(Q)}{l}\Big)\eps.
  \end{align*}
  This bound implies the continuity.  The case in which $N>2$ follows
  because, moving all points by at most $\eps$, the change in all the
  regions is upper bounded by $\bigcup_{1\leq i ,j\le N }Y_{i j}(\eps)$,
  which vanishes as $\eps\to0^+$.

  The last remaining step is to prove claims~\ref{item:TCont}
  and~\ref{item:TepsCont}. We focus on claim~\ref{item:TepsCont} and show
  claim~\ref{item:TCont} as a byproduct.  Let $v\in\part$ and let
  $\widehat{v}=\Tdelta_{ij}(v)$. According to the definitions in
  Subsection~\ref{sec:SmoothedMap}, $\widehat{v}$ is characterized by the
  sets $\widehat{R}_i$ and $\widehat{R}_j$.  Recall that these two sets
  depend on the sets $R_i$ $R_j$, on the scalar $\beta_{ij}(v)$ and on the
  points $\widehat{p}_i$ and $\widehat{p}_j$.  One can see that
  $\beta_{ij}(v)$ is a continuous function of its arguments $v_i$ and
  $v_j$.  Hence, it suffices to show that also $\widehat{R}_i$,
  $\widehat{R}_j$ and $\widehat{p}_i$, $\widehat{p}_j$ depend continuously
  on $v_i$ and $v_j$.  To do this, introduce $v'\in\part$ and compute
  $\widehat{v}'=\Tdelta_{ij}(v')$.  Assume $\Cd(v_i)\neq \Cd(v_j)$ and
  $\Cd(v_i')\neq \Cd(v_j')$.  Analogously to how we defined $R_i$ and
  $R_j$, we now define the regions $R'_i= v'_i \intersect
  \halfs(\Cd(v'_j),\Cd(v'_i))$ and $R'_j= v'_j \intersect
  \halfs(\Cd(v'_i),\Cd(v'_j))$.  We aim to upper bound the composite
  distance $d_\ominus(R_i, R'_i)+d_\ominus(R_j, R'_j)$.
  Observe that this composite distance depends on the difference of the two
  argument regions $v_i, v'_i$ and $v_j, v'_j$ both directly and indirectly
  via the induced difference between the centroids.  Recalling the proof of
  claim~\ref{item:V-cont}, let $\eps$ be an upper bound on the displacement
  between the two centroids. Then the region $Y(\eps)=\setdef{x\in Q}{
    |\|x-\Cd(v_i)\| - \|x-\Cd(v_j)\| \le 2\eps}$ needs to be included in
  the upper bound on the composite distance. Combining these considerations
  we obtain
    \begin{multline}\label{eq:Up-Bound}
      d_\ominus(R_i, R'_i)+d_\ominus(R_j, R'_j)\le
      \left(d_\ominus(v_i, v'_i)+d_\ominus(v_j, v'_j) \right) \\ +
      \measure\big( Y(\max\{\|\Cd(v_i)-\Cd(v'_i)\|,
        \|\Cd(v_j)-\Cd(v'_j)\| \}) \big).
    \end{multline}
    Clearly, if $d_\ominus(v_i,v'_i)\to 0$ and $d_\ominus(v_j,v'_j)\to 0$,
    then $\|\Cd(v_i)-\Cd(v'_i)\|\to 0$ and $\|\Cd(v_i)-\Cd(v'_i)\|\to 0$
    and, in turn, also $d_\ominus(R_i, R'_i)+d_\ominus(R_j, R'_j)\to
    0$. Hence, we can argue that the sets $R_i$ and $R_j$ depend
    continuously on the regions $v_i$ and $v_j$.  This is enough to prove
    statement~\ref{item:TCont}, provided the two regions $v_i$ and $v_j$
    have distinct centroids.  Moreover a direct consequence of this fact is
    that also the points $\widehat{p}_i$ and $\widehat{p}_j$ depend
    continuously on $v_i$ and $v_j$. Finally, observe that in the limit
    case $\Cd(v_i)=\Cd(v_j)$ the continuity of $R_i$ and $R_j$ is captured
    by the fact that $\beta_{ij}(v)$ is a continuous function of
    $\Cd(v_i),\Cd(v_j)$ and that $\beta_{ij}(v)=0$ if $\Cd(v_i)=\Cd(v_j)$.
    The continuity of $R_i$, $R_j$, $\beta_{ij}(v)$, $\widehat{p}_i$ and
    $\widehat{p}_j$ imply also the continuity of $\widehat{R}_i$ and
    $\widehat{R}_j$ and, in turn, of $\Tdelta_{ij}$.
\end{proof}

\subsection{Convergence proofs}
\label{sec:ConvergProofs}
In view of the identification between $N$-partitions and their equivalence
classes introduced at the beginning of this section, we are now ready to
complete the proof of the convergence results presented in
Section~\ref{sec:AlgoConverg}.

We start by clarifying the precise meaning of convergence in
Theorems~\ref{th:T-persistent-gossip}
and~\ref{th:Tdelta-persistent-gossip}. Specifically, we say that a sequence
of partitions $\seqdef{v(t)}{t\in\integernonnegative} \subset \part$
converges to a set of partitions $X\subset\part$ if the symmetric distance
from $\seqdef{v(t)}{t\in\integernonnegative}$ to $X$ converges to zero,
that is,
\begin{equation*}
  \lim_{t\to\infty} \enspace \inf \setdef{d_\ominus(v(t),x)} {x\in{X}} =  0.
\end{equation*}
\begin{proof}[Proof of Theorem~\ref{th:Tdelta-persistent-gossip}]
  We prove the deterministic statement~(i).  We start by observing that,
  through the canonical projection, the evolution
  $\seqdef{v(t)}{t\in\integernonnegative} \subset \part$ of
  $\map{\Tdelta}{\part}{\part}$ can be mapped into the evolution
  $\seqdef{v^*(t)}{t\in\integernonnegative} \subset \part^*$ of
  $\map{\Tdelta}{\part^*}{\part^*}$.  We aim to apply
  Theorem~\ref{th:StrongerLasalle} to the dynamical system
  $\map{\Tdelta}{\part^*}{\part^*}$ and its evolution
  $\seqdef{v^*(t)}{t\in\integernonnegative} \subset \part^*$.  In what
  follows, our goal is to verify whether
  Assumptions~\ref{item:Compactness},~\ref{item:Lyapunov},~\ref{item:continuousIterate}
  and~\ref{item:persistentChoice} of Theorem~\ref{th:StrongerLasalle} are
  satisfied.

  Since $\seqdef{v(t)}{t\in\integernonnegative}$ is non-vanishing and
  finitely convex by assumption, it follows that there exists
  $\ell\in\natural$ such that the $\Omega$-limit set of
  $\seqdef{v(t)}{t\in\integernonnegative}$ is contained in $\CL\cap\part$,
  that is, $\Omega(v(t))\subseteq \CL\cap \part$. This implies also that
  $\Omega(v^*(t))\subseteq \CL^*\cap\part^*$.  As stated in
  Theorem~\ref{th:CL}, $\CL^*$ is compact in the topology induced by the
  metric $d_\ominus$. Hence, even though $\CL^*$ is not strongly positive
  invariant for $\Tdelta$, the weaker version of
  Assumption~\ref{item:Compactness} of Theorem~\ref{th:StrongerLasalle}, as
  given in Remark~\ref{rem:NoStrongly}, holds true for the sequence
  $\seqdef{v^*(t)}{t\in\integernonnegative} \in \part^*$.  Now, as one can
  deduce from Theorem~\ref{th:EveryThingCont}\ref{item:Hcentroid-is-cont}
  and Lemma~\ref{lem:HcentroidPairwiseDecr}, the function $\Hcentroid$
  satisfies the Assumption~\ref{item:Lyapunov} of
  Theorem~\ref{th:StrongerLasalle}, thus playing the role of a Lyapunov
  function for the dynamical system $\Tdelta$. Moreover, from
  Theorem~\ref{th:EveryThingCont}\ref{item:TepsCont}, note that the system
  evolves through maps that are continuous in $\part^*$ with respect to the
  metric $d_\ominus$: thus the Assumption~\ref{item:continuousIterate} of
  Theorem~\ref{th:StrongerLasalle} is satisfied.  Finally observe that the
  Assumption~\ref{item:persistentChoice} of
  Theorem~\ref{th:StrongerLasalle} corresponds to the assumption of uniform
  persistency in Theorem~\ref{th:Tdelta-persistent-gossip}.  Therefore, we
  conclude that the evolution $\seqdef{v^*(t)}{t\in\integernonnegative}$
  converges to the intersection of the fixed points of the maps
  $\Tdelta_{ij}$, for all $i, j \in \until{N}$, $j\neq i$. According to
  Lemma~\ref{lem:fixed-points-cvp}, this intersection coincides with the
  set of mixed centroidal Voronoi partitions up to sets of measure zero.

  The proof of the stochastic statement~(ii) follows the same lines,
  applying Theorem~\ref{th:StrongerLasalleProb} instead of
  Theorem~\ref{th:StrongerLasalle}.
\end{proof}

\begin{proof}[Proof of Theorem~\ref{th:T-persistent-gossip}]
  The proof of this result follows the lines of the proof above with two
  distinctions. The distinctions come from the assumption that the
  evolution $v$, in addition to being non-vanishing and finitely convex, is
  also distinct centroidal.  The first distinction is as follows. In order
  to apply the Krasovskii-LaSalle invariance principle we require the
  continuity property stated in
  Theorem~\ref{th:EveryThingCont}\ref{item:TCont}. Additionally, we note
  that the space of finitely convex and distinct centroidal partitions
  \begin{equation*}
    \setdef{v\in\part\intersect \CL^N}{
      \|\Cd(v_i)-\Cd(v_j)\|\geq \eps \text{ for all }  i\neq j}
  \end{equation*}
  is a closed and hence compact subspace of $\part$.  The second
  distinction is as follows.  Since we rule out the case of coincident
  centroids, we can infer convergence to centroidal Voronoi partitions
  instead of convergence to mixed centroidal Voronoi partitions; see
  Lemma~\ref{lem:fixed-points-cvp}\ref{item:from-mixed-to-cvt}.
\end{proof}

\section{Conclusions}
\label{sec:Outro}
In summary, we have introduced novel coverage, deployment and partitioning
algorithms for robotic networks with minimal communication requirements. To
analyze our proposed algorithms, we have developed and characterized (1)
intuitive versions of the Krasovskii-LaSalle Invariance Principle for
deterministic and stochastic switching systems, (2) relevant topological
properties of the space of partitions, and (3) useful continuity properties
of a number of geometric and multicenter functions.

We believe there remain interesting open issues in the study of gossiping
robots and of dynamical systems on the space of partitions.  We are keen on
extending these ideas to non-convex complex environments and discrete
environments such as graphs. Following Remark~\ref{rem:disk-covering} we
plan to study gossip coverage algorithms for more general multicenter
functions, including nonsmooth, anisotropic and inhomogeneous functions.
Additionally, we plan to investigate gossip coverage algorithms capable of
adapting to time-varying scenarios such as problems in which robotic agents
arrive to and depart from the network.
Finally, inspired by stigmergy in territorial animals, we plan to design
communication protocols for multiagent systems based on the ability of
leaving messages in the environment.

\appendix

\section{A robotic network implementation of gossip coverage algorithms}
\label{app:robotic-implementation}
In Subsection~\ref{subsec:distribcoverage+delaunaynetwork} we discussed
coverage control algorithm for groups of robots with synchronized and
reliable communication along all edges of a Delaunay graph; then in
Subsection~\ref{subsec:gossip-coverage-algo} we introduced our gossip
coverage algorithms. Here we discuss in full detail one possible way of
implementing the partition-based gossip coverage algorithm in a robotic
network with weak communication requirements.

We consider a group of agents all having the following capabilities: (C1)
each agent $i\in\until{N}$ knows its position and moves at positive speed
$u_i$ to any position in the compact convex environment $Q\subset\real^2$;
(C2) each agent may store an arbitrary number of locations in $Q$ and has a
clock that is not necessarily synchronized with other agents' clocks
(specifically, we assume same clock skew, but different clock offsets among
the agents); and (C3) if any two agents are within distance $\rcomm$ of
each other for some positive duration of time, then they exchange
information at the sample times of a Poisson process with intensity
$\lcomm$.

The \emph{random destination \& wait algorithm} is described as follows.
Given a parameter $\eps<\rcomm/4$, each agent $i\in\until{N}$ maintains in
memory a dominance region $v_i$ and determines its motion by repeatedly
performing the following three actions over periods of time that we label
\emph{epochs}:
\begin{algorithmic}[1]
  \STATE instantaneously agent $i$ selects uniformly randomly a
  \emph{destination point} $q_i$ in the set $\setdef{q\in\real^2}
  {\dist\big(q,\partial{v_i}\setminus\partial{Q}\big)\leq\eps}$;

  \STATE agent $i$ moves in such a way as to reach point $q_i$ in time
  precisely equal to $\duration=\diam(Q)/\min\{u_1,\dots,u_N\}$; and

  \STATE agent $i$ waits at point $q_i$ for a duration that equals either
  $\duration$ with probability $1/2$ or $2\duration$ with probability
  $1/2$.
\end{algorithmic}
We have assumed that agents may reach locations at a distance up to $\eps$
away from $Q$; this assumption can be removed at the cost of additional
notation.

The random destination \& wait algorithm is meant to be implemented
concurrently with the modified gossip coverage algorithm -- where the
parameter $\delta$ is selected to satisfy $\delta<\rcomm/4$. The two
algorithms jointly determine the evolution of the agents positions and
dominance regions as follows.  If any two agents $i$ and $j$ are within
communication range $\rcomm$ at any instant of time and for some positive
duration, then, at each sample time of the corresponding communication
Poisson process, the two agents exchange sufficient information to update
their respective regions $v_i$ and $v_j$ via the modified gossip coverage
map $\Tdelta_{ij}$.

\begin{proposition}[Random destination \& wait ensures persistent random
 gossip]
 \label{prop:I-promise-this-is-the-last-one}
 Consider a group of agents with capacities (C1), (C2) and (C3) and
 parameters $u_i$, $\rcomm$ and $\lcomm$. Assume the agents implement the
 random destination \& wait algorithm and the modified gossip coverage
 algorithm with parameter $\eps<\rcomm/4$ and $\delta<\rcomm/4$.  The
 following statements hold:
 \begin{enumerate}
 \item the sequence of applications of the modified gossip coverage map is
   a randomly persistent stochastic process; and
 \item the resulting evolution $\map{v}{\realnonnegative}{\part}$,
   conditioned upon being non-vanishing and finitely convex, converges
   almost surely to the set of mixed centroidal Voronoi partitions.
 \end{enumerate}
\end{proposition}

\begin{proof}
  Assume that at time $s\in\realnonnegative$, the regions $v_i$ and $v_j$
  are at distance $\delta$ or less. Pick any two points
  $p_i\in\partial{v_i}\setminus\partial{Q}$ and
  $p_j\in\partial{v_j}\setminus\partial{Q}$ such that $\|p_i-p_j\|\leq 2
  \delta$ and define two balls of radius $\eps$ centered at $p_i$ and $p_j$
  by $B_i=\setdef{q\in\real^2}{\|q-p_i\|\leq\eps}$ and
  $B_j=\setdef{q\in\real^2}{\|q-p_j\|\leq\eps}$.  By construction, each
  point in $B_i$ is at most at distance $2\eps+2\delta<\rcomm$ from each
  point in $B_j$.
  We claim that, independently of the agent positions and the environment
  partitions at time $s$ and of their evolution throughout the past time
  interval $[0,s]$, there exists a duration $\Delta\in\realnonnegative$
  such that the meeting event ``agent $i$ is in $B_i$ and agent $j$ is in
  $B_j$ during the time interval $[s+\Delta,s+\Delta+d]$'' has probability
  that is lower bounded by a positive constant. If these meeting events have
  positive probability, then communication events have positive probability
  and therefore the proposition follows.

  We prove the claim as follows.  At time $s$, we say that an agent is in
  \OFF state if it is at the initial time of an epoch or it is moving to
  its next destination point.  We say it is in $\ON_1$ state if it has just
  reached a destination point or it has been waiting at the destination
  point for a duration of time less than $\duration$.  We say it is in
  $\ON_2$ state if it has been waiting at the destination point for a
  duration of time greater than or equal to $\duration$ and less than
  $2\duration$.  Note that each agent remains in each one of the three
  states $\{\OFF,\ON_1,\ON_2\}$ for a duration of time $\duration$.  The
  transitions between states are regulated by a Markov chain with three
  states $\{\OFF,\ON_1,\ON_2\}$ and with (column-stochastic) transition
  matrix
  \begin{equation*}
    A =
    \begin{bmatrix}
      0 & 1/2 & 1 \\
      1 & 0 & 0 \\
      0 & 1/2 & 0
    \end{bmatrix}.
  \end{equation*}
  For now, suppose that at time $s$, agent $i$ is at the initial time of an
  interval of duration $\duration$ in either state \OFF, or state $\ON_1$
  or state $\ON_2$.  From the matrices
  \begin{equation*}
    A^3 =
    \begin{bmatrix}
      1/2 & 1/4 & 1/2 \\
      1/2 & 1/2 & 0 \\
      0 & 1/4 & 1/2
    \end{bmatrix}, \quad \text{and} \quad
    A^5 =
    \begin{bmatrix}
      1/2 & 3/8 & 1/4 \\
      1/4 & 1/2 & 1/2 \\
      1/4 & 1/8 & 1/4
    \end{bmatrix},
  \end{equation*}
  we establish the following three facts:
  \begin{enumerate}
  \item if agent $i$ is at state $\OFF$ at time $s$, then at time
    $s+3\duration$ it is at state $\OFF$ with probability $1/2$ and at time
    $s+5\duration$ it is at state $\ON_2$ with probability $1/4$;
  \item if agent $i$ is at state $\ON_1$ at time $s$, then at time
    $s+3\duration$ it is at state $\OFF$ with probability $1/4$ and at time
    $s+5\duration$ it is at state $\ON_2$ with probability $1/8$; and
  \item if agent $i$ is at state $\ON_2$ at time $s$, then at time
    $s+3\duration$ it is at state $\OFF$ with probability $1/2$ and at time
    $s+5\duration$ it is at state $\ON_2$ with probability $1/4$.
  \end{enumerate}
  We now return to the two agents $i$ and $j$.  According to the robot
  capability (C2), at time $s$ the clocks of agent $i$ and agent $j$,
  denoted by $t_i(s)$ and $t_j(s)$ respectively, are distinct and can be
  written as an integer number of durations $d$ plus a quotient;
  specifically $t_i(s)=k_i\duration +\eta_i$, and $t_j(s)= k_j\duration
  +\eta_j$, with $k_i, k_j \in \integernonnegative$, $\eta_1, \eta_2 \in
  [0, \duration[$, and possibly $k_i\neq k_j$ and $\eta_i \neq
  \eta_j$. Without loss of generality, assume $\eta_i\leq \eta_j$. Facts
  (i), (ii), and (iii) and the assumption that clock skew is equal among
  agents, jointly imply that, with probability greater or equal to
  $(1/8)^2=1/64$, agents $i$ and $j$ during the interval
  $\left[s+4\duration-\eta_i, s+5\duration-\eta_i\right]$ are waiting at
  some destination point chosen, respectively, at time
  $s+3\duration-\eta_i$ and $s+3\duration-\eta_j$.  If agent $i$ is in
  $B_i$ and agent $j$ is inside $B_j$ during the interval
  $\left[s+4\duration-\eta_i, s+5\duration-\eta_i\right]$, then they
  communicate at least once with probability $\pcomm :=
  1-\ee^{-\duration\lcomm}>0$. Let us now find a lower bound on the
  probability of agents $i$ being in $B_i$ at time $s+4\duration-\eta_i$;
  clearly this happens if the following two facts occur:
  \begin{itemize}
  \item during the interval $\left[s, s+3 \duration-\eta_i\right]$
    agent $i$ does not communicate with any other agent and hence its
    dominance region $v_i$ remains unchanged; and
  \item at time instant $s+3 \duration-\eta_i$ agent $i$ is at the
    beginning of a new epoch (state \OFF) and selects a destination point
    in $B_i$.
  \end{itemize}
  Since during the interval $\left[s-\eta_i, s+3 \duration-\eta_i\right]$
  agent $i$ spends a time duration of at most $2\duration$ at distance less
  than~$\delta$ from any other agents, then the probability that agent $i$
  does not communicate in that interval of time is lower bounded by
  $(1-\pcomm)^2=\ee^{-2d\lcomm}$.  Moreover, provided that at time instant
  $s+3\duration-\eta_i$ agent $i$ is at the beginning of a new epoch (state
  $\OFF$) and that the $v_i$ is unchanged until time instant $s$, we have
  that the destination point for agent $i$ is selected uniformly randomly
  so that it belongs to $B_i$ with probability
  $\operatorname{area}(B_i)/\operatorname{area}\big(\setdef{q\in\real^2}
  {\dist\big(q,\partial{v_i}\setminus\partial{Q}\big)\leq\eps}\big) \geq
  \parea :=
  \pi\eps^2/\big(\operatorname{area}(Q)+\eps\operatorname{perimeter}(Q)\big)$.
  Similar considerations hold also for agent~$j$.

  In summary, if agents $i$ and $j$ are neighbors at time $s$, then they
  communicate at least once during the interval $\left[s+4\duration-\eta_i,
    s+5\duration-\eta_i\right]$ with probability lower bounded by
  $(1/64) \parea^2 (1-\pcomm)^2 \pcomm>0$.  If agents $i$ and $j$ are at
  distance greater than $\delta$ at time $s$ (and therefore they are not
  neighbors), then the modified gossip coverage map is the identity map so
  that their communication is immaterial.
\end{proof}

\section{A counterexample showing the necessity of uniformly persistent
  switches}
\label{appA:counterexample}
Theorem~\ref{th:StrongerLasalle} contains a persistent switching
conditions, that is, it requires the existence of $D\in \natural$ such that
every map $T_i$, $i \in \until{m}$, is applied at least once within every
interval $[n, n +D[$, $n \in \integernonnegative$. This appendix contains
an example proving the necessity of this condition.

Consider the plane in polar coordinates $X=\realpositive \times {[0,2\pi[}
\union \{0,0\}$.  Define the standard metric
$\map{d}{X\times{X}}{\realnonnegative}$ as follows: let $(\rho_1,
\theta_1)$, $(\rho_2, \theta_2)$ be any pair of elements of $X$ and
  \begin{equation*}
    d((\rho_1, \theta_1), (\rho_2,
    \theta_2))
    =\sqrt{\left(\rho_1\cos \theta_1-\rho_2 \cos
        \theta_2\right)^2+ \left(\rho_1\sin \theta_1-\rho_2 \sin
        \theta_2\right)^2}.
  \end{equation*}
  Consider now the continuous maps $\map{T_i}{X}{X}$, $i\in\{1,2\}$,
  defined by respectively
  \begin{align*}
    T_1(\rho,\theta)&=
    \begin{cases}
      (\rho^2, \theta), & \text{if }  0 \leq \rho \leq 1,\\
      \left(\frac{2\rho-1}{\rho},\, (\theta+\rho-1)\mod 2\pi\right), &
      \text{if } \rho>1 ,
    \end{cases}
    \\
    T_2(\rho,\theta)&=
    \begin{cases}
      \left( (1-\sin\theta) \rho,\theta\right), \qquad\qquad\qquad &\text{if } 0\leq \theta
      \leq \pi, \\
      \left(\rho, \theta\right) & \text{if } \pi \leq \theta \leq 2\pi.
    \end{cases}
  \end{align*}
  Define $\setmap{T}{X}{X}$ by $T(\rho,\theta) = \left\{T_1\left(\rho,
      \theta\right), T_2\left(\rho,\theta\right)\right\}$ and the function
  $\map{U}{X}{\realnonnegative}$ by
  $U\left(\rho,\theta\right)=\rho$. Observe that $U$ is continuous and
  non-increasing along $T$. Assume now that there exists $D\in\natural$ such
  that, for any $n\in\integernonnegative$, there exist $n_1$ and $n_2$
  within the interval $]n, n+D]$ such that $x_{n_1+1}=T_1(x_{n_1})$ and
  $x_{n_2+1}=T_2(x_{n_2})$. Then, by Theorem \ref{th:StrongerLasalle}, the
  $\omega$-limit set of each evolution of $T$ is a subset of
  \begin{equation}
    \label{eq:omegaset-with-uniformity}
    \setdef{(\rho, \theta)\in X} {\rho=1, \pi \leq \theta \leq 2 \pi} \cup
    \{0,0\}.
  \end{equation}

  Next, we relax the condition that the map $T_2$ is applied at least once
  inside each interval of arbitrary amplitude $D$ and we show that there
  exists one sequence that does not converge to the $\omega$-limit set in
  equation~\eqref{eq:omegaset-with-uniformity}. To this aim, assume the
  sequence $\seqdef{ (\rho(n), \theta(n))}{n\in\integernonnegative}$
  satisfies
  \begin{enumerate}
  \item $\rho(0)>1$;
  \item  $(\rho(1), \theta(1))= T_1\left(\rho(0), \theta(0)\right)$ and
  \item $(\rho(n+1), \theta(n+1))= T_2\left(\rho(n), \theta(n)\right)$ if
    and only if $\pi \leq \theta(n) \leq 2 \pi$ and $(\rho(n), \theta(n))=
    T_1\left(\rho(n-1), \theta(n-1)\right)$.
  \end{enumerate}
  Note that if $\pi \leq \theta(n) \leq 2 \pi$, then $ T_2\left(\rho(n),
    \theta(n)\right)=(\rho(n), \theta(n))$. Therefore, the evolution
  $\left\{(\rho(n), \theta(n))\right\}$ equals $\left\{(\widehat{\rho}(n),
    \widehat{\theta}(n))\right\}$ where $(\widehat{\rho}(0),
  \widehat{\theta}(0))=(\rho(0),\theta(0))$ and $(\widehat{\rho}(n),
  \widehat{\theta}(n))=T_1^n\left(\widehat{\rho}(0),
    \widehat{\theta}(0)\right)$. Regarding this new sequence, observe that
  \begin{gather}
    1<\widehat{\rho}(i)<2 \text{ and } \widehat{\rho}(i+1)<
    \widehat{\rho}(i), \text{ for all }i \geq 1, \label{[(a)]}
    \\
    0<\widehat{\theta}(i+1)-\widehat{\theta}(i)<\pi, \text{ for all }i\geq
    1, \text{ and } \lim_{i\to \infty}
    \left(\widehat{\theta}(i+1)-\widehat{\theta}(i)\right)=0, \label{[(b)]}
    \\
    \lim_{r\to \infty} \sum_{i=1}^r
    \left(\widehat{\theta}(i+1)-\widehat{\theta}(i)\right)= \lim_{r\to
      \infty} \sum_{i=1}^r \left(\widehat{\rho}(i)-1\right)= \lim_{r\to
      \infty} \sum_{i=1}^r
    \left(\frac{1}{\widehat{\rho}(1)-1}+i-1\right)^{-1}
    \!\!=\infty, \label{[(c)]}
  \end{gather}
  where the equality $\widehat{\rho}(i)-1=
  (\frac{1}{\widehat{\rho}(1)-1}+i-1)^{-1}$ can be proved by induction over
  $i$.  Properties~\eqref{[(a)]}, \eqref{[(b)]}, and~\eqref{[(c)]} ensure
  that there exists a sequence $\setdef{n_h}{h\in\integernonnegative}$ such
  that $(\rho(n_h), \theta(n_h))=T_2\left(\rho(n_h-1),
    \theta(n_h-1))\right)$ for all $h\in \integernonnegative$, and
  $(\rho(t), \theta(t))=T_1\left(\rho(t-1), \theta(t-1))\right)$ if $t
  \notin \setdef{n_h}{h\in\integernonnegative}$. Moreover, we have that
  $\lim_{h\to \infty}(n_{h+1}-n_h)=\infty$. In other words, both the maps
  $T_1$ and $T_2$ are applied \emph{infinitely often} along the evolution
  described by $\left\{(\rho(n), \theta(n))\right\}$, but there does not
  exists $D\in \natural$ such that $T_2$ is applied at least once within
  each interval $[n, n+D]$, $n\in\integernonnegative$.  Observe that, in
  this case, $(\rho(n), \theta(n))$ converges to the set
  $\setdef{(1,\theta)}{\theta\in[0,2\pi]}\subset{X}$.  This set is
  different from the $\omega$-limit set in
  equation~\eqref{eq:omegaset-with-uniformity}.

\section{Discontinuity of the multicenter function in the Hausdorff metric}
\label{app:discontinuous}
To see that $\map{\Hcentroid}{\part^*}{\realnonnegative}$ is not
Hausdorff-continuous, consider the sequence of $2$-partitions
$\seqdef{v(t)}{t\in\integernonnegative}$ of the interval $[-1, 1] \subseteq
\real$ defined by
\begin{equation*}
  v_1(t)=\left[-1, -1+\frac{1}{2^{t+1}}\right]\cup
  \bigcup_{h=-(2^{t-1}-1)}^{2^{t-1}-1}\left[\frac{h}{2^{t-1}}-\frac{1}{2^{t+1}},
    \frac{h}{2^{t-1}}+\frac{1}{2^{t+1}}\right]\cup
  \left[1-\frac{1}{2^{t+1}},1\right],
\end{equation*}
and by $v_2(t)=\closure{[-1,1]\setminus v_1(t)}$.  Note that both sequences
$\seqdef{v_1(t)}{t\in\integernonnegative}$ and
$\seqdef{v_2(t)}{t\in\integernonnegative}$ converge to $[-1,1]$, and that
$\Cd(v_1(t))=\Cd(v_2(t))=\Cd([-1,1])=0$, for all $t\ge 0$. Hence, for
$\phi(s)=1$ and $f(s)=s$, we compute
\begin{equation*}
  \Hone(0,v_1(t))=\int_{v_1(t)}|x|dx=2 \int_{v_1(t)\intersect [0,1]}|x|dx= 2
  \frac{1}{2}=1,
\end{equation*}
and consequently $\Hcentroid(v(t))=
2$, while $\Hcentroid(\lim_{t\to\infty}v(t))= 2\Hone(0,[-1,1])=4$. This
shows that $\lim_{t \to \infty}\Hcentroid(v(t))\neq \Hcentroid(\lim_{t\to
  \infty} v(t))$.

\bibliographystyle{siam}

\end{document}